\documentclass[12pt]{amsart}
\usepackage{amssymb,latexsym}
\usepackage{array}
\usepackage{tikz-cd}
\usepackage{amsmath, amssymb, amsfonts, amstext, amsthm, amscd}
\usepackage[mathscr]{euscript}
\usepackage{enumerate}
\usepackage{tikz,color,soul}
\usepackage[Symbolsmallscale]{upgreek}
\usetikzlibrary{arrows.meta, positioning,arrows}
\usetikzlibrary{matrix,arrows}
\usepackage[english]{babel}
\usepackage{chngcntr}
\usepackage{textgreek}

\begin{document}

\thispagestyle{empty} 

\title{Equivariant $K$-theory of cellular toroidal embeddings}

%\author{Alexis Tchoudjem and V. Uma}

\author[A. Tchoudjem]{Alexis Tchoudjem}
\address{Institut Camille Jordan
Université Lyon 1
43, boulevard du ONZE NOVEMBRE
69622 Villeurbanne Cedex
FRANCE}
\email{tchoudjem@math.univ-lyon1.fr}

\author[V. Uma]{Vikraman Uma}
\address{Department of Mathematics, Indian Institute of Technology, Madras, Chennai 600036, India
}
\email{vuma@iitm.ac.in}


%\subjclass{55N15, 14M15, 19L99}


\date{}

\begin{abstract}
  In this article we describe the
  $G_{comp}\times G_{comp}$-equivariant topological $K$-ring of a {\em
    cellular} toroidal embedding $\mathbb{X}$ of a complex connected
  reductive algebraic group $G$. In particular, our results extend the
  results in \cite{u1} and \cite{u2} on the regular embeddings of $G$,
  to the equivariant topological $K$-ring of a larger class of
  (possibly singular) cellular toroidal embeddings. They are also a
  topological analogue of the results in \cite{gon} on the operational
  equivariant algebraic $K$-ring, for cellular toroidal
  embeddings.\end{abstract}

\maketitle

\def\theequation
  {\arabic{section}.\arabic{equation}}

\newcommand{\codim}{\mbox{{\rm codim}$\,$}}
\newcommand{\stab}{\mbox{{\rm stab}$\,$}}
\newcommand{\lr}{\mbox{$\rightarrow$}}

\newcommand{\be}{\begin{equation}}
\newcommand{\ee}{\end{equation}}

\newtheorem{guess}{Theorem}[section]
\newcommand{\bth}{\begin{guess}$\!\!\!${\bf }~}
\newcommand{\eeth}{\end{guess}}
\renewcommand{\bar}{\overline}
\newtheorem{propo}[guess]{Proposition}
\newcommand{\bpropo}{\begin{propo}$\!\!\!${\bf }~}
\newcommand{\epropo}{\end{propo}}

\newtheorem{assum}[guess]{Assumption}
\newcommand{\bas}{\begin{assum}$\!\!\!${\bf }~}
\newcommand{\eas}{\end{assum}}

\newtheorem{lema}[guess]{Lemma}
\newcommand{\blem}{\begin{lema}$\!\!\!${\bf }~}
\newcommand{\elem}{\end{lema}}

\newtheorem{defe}[guess]{Definition}
\newcommand{\bdefe}{\begin{defe}$\!\!\!${\bf }~}
\newcommand{\edefe}{\end{defe}}

\newtheorem{coro}[guess]{Corollary}
\newcommand{\bcor}{\begin{coro}$\!\!\!${\bf }~}
\newcommand{\ecor}{\end{coro}}

\newtheorem{rema}[guess]{Remark}
\newcommand{\brem}{\begin{rema}$\!\!\!${\bf }~\rm}
\newcommand{\erem}{\end{rema}}

\newtheorem{exam}[guess]{Example}
\newcommand{\beg}{\begin{exam}$\!\!\!${\bf }~\rm}
\newcommand{\eeg}{\end{exam}}

\newtheorem{notn}[guess]{Notation}
\newcommand{\bnot}{\begin{notn}$\!\!\!${\bf }~\rm}
\newcommand{\enot}{\end{notn}}

\newcommand{\ch}{{\mathcal H}}
\newcommand{\cf}{{\mathcal F}}
\newcommand{\cd}{{\mathcal D}}
\newcommand{\cR}{{\mathcal R}}
\newcommand{\cv}{{\mathcal V}}
\newcommand{\cn}{{\mathcal N}}
\newcommand{\lra}{\rightarrow}
\newcommand{\ra}{\rightarrow}
\newcommand{\blr}{\Big \rightarrow}
\newcommand{\da}{\Big \downarrow}
\newcommand{\ua}{\Big \uparrow}
\newcommand{\hra}{\mbox{{$\hookrightarrow$}}}
\newcommand{\rt}{\mbox{\Large{$\rightarrowtail$}}}
\newcommand{\dua}{\begin{array}[t]{c}
\Big\uparrow \\ [-4mm]
\scriptscriptstyle \wedge \end{array}}
\newcommand{\ctext}[1]{\makebox(0,0){#1}}
\setlength{\unitlength}{0.1mm}
\newcommand{\cl}{{\mathcal L}}
\newcommand{\cp}{{\mathcal P}}
\newcommand{\ci}{{\mathcal I}}
\newcommand{\bz}{\mathbb{Z}}
\newcommand{\cs}{{\mathcal s}}
\newcommand{\ce}{{\mathcal E}}
\newcommand{\ck}{{\mathcal K}}
\newcommand{\cz}{{\mathcal Z}}
\newcommand{\cg}{{\mathcal G}}
\newcommand{\cj}{{\mathcal J}}
\newcommand{\cc}{{\mathcal C}}
\newcommand{\ca}{{\mathcal A}}
\newcommand{\cb}{{\mathcal B}}
\newcommand{\cx}{{\mathcal X}}
\newcommand{\co}{{\mathcal O}}
\newcommand{\bq}{\mathbb{Q}}
\newcommand{\bt}{\mathbb{T}}
\newcommand{\bh}{\mathbb{H}}
\newcommand{\br}{\mathbb{R}}
\newcommand{\bl}{\mathbf{L}}
\newcommand{\wt}{\widetilde}
\newcommand{\im}{{\rm Im}\,}
\newcommand{\bc}{\mathbb{C}}
\newcommand{\bp}{\mathbb{P}}
\newcommand{\ba}{\mathbb{A}}
\newcommand{\spin}{{\rm Spin}\,}
\newcommand{\ds}{\displaystyle}
\newcommand{\tor}{{\rm Tor}\,}
\newcommand{\bff}{{\bf F}}
\newcommand{\bs}{\mathbb{S}}
\def\ns{\mathop{\lr}}
\def\nssup{\mathop{\lr\,sup}}
\def\nsinf{\mathop{\lr\,inf}}
\newcommand{\tT}{{\widetilde{T}}}
\newcommand{\tG}{{\widetilde{G}}}
\newcommand{\tB}{{\widetilde{B}}}
\newcommand{\tC}{{\widetilde{C}}}
\newcommand{\tW}{{\widetilde{W}}}
\newcommand{\tphi}{{\widetilde{\Phi}}}

%%%%%%%%%%

\newcommand{\ie}{\emph{i.e.}}.
%%%%%%%%%%%%%

\noindent

\section{Introduction}\label{Introduction}

Let $G$ denote a complex connected reductive algebraic group.  Let $C$ denote the center of $G$ and let $G_{ad}=G/C$ be the corresponding 
semisimple adjoint group. 

%Let $\bar{G_{ad}}$ denote the {\it wonderful
 % compactification} of $G_{ad}$ (see \cite{dp}).

An {\it equivariant embedding} $\mathbb{X}$ of $G$ is a normal
$G\times G$-variety containing a point $x$ such that the orbit
$(G\times G)\cdot x$ is open and isomorphic to $G\cong G\times
1$. When $G=T$ is a torus, the $T$-embeddings are nothing but the $T$-toric
varieties.  We say that the equivariant embedding $\mathbb{X}$ of $G$
is {\em toroidal} if the quotient map $p:G\lra G_{ad}$ extends to a
unique $G\times G$-equivariant morphism
$p:\mathbb{X}\lra \bar{G_{ad}}$, where $\bar{G_{ad}}$ denotes the
wonderful compactification of $G_{ad}$ defined by De Concini and
Procesi (see \cite{dp}). Here the $G\times G$-action on $\bar{G_{ad}}$ is
via its projection to $G_{ad}\times G_{ad}$.

The smooth toroidal embeddings are exactly the regular embeddings (see
\cite[Chapter 5, Theorem 29.2]{t}, \cite{bdp}, \cite[Section
1.4]{Br1}).

Let $\mathbb{X}$ denote a toroidal embedding of $G$.  Let $T\subseteq G$ a maximal torus
of $G$. Let $W$ denote the Weyl group of $(G,T)$. Let
$X=\overline{(T\times T)\cdot x}$ denote the closure of
$T\cong T\times 1$ in $\mathbb{X}$.  The fan $\cf$ associated to $X$
is a subdivision of the Weyl chambers in $X_*(T)_{\mathbb{R}}$. Let
$T\subseteq B\subseteq G$ be a Borel subgroup
of $G$. Let
$X^+$ denote the toric variety whose fan $\cf^+$ is the subfan of
$\cf$ consisting of the cones in the positive Weyl chamber. In
particular, $\cf^+$ is a subdivision of the positive Weyl chamber
$\mathcal{C}^+$ associated to $(G\supseteq B\supseteq T)$ (see
Section \ref{celltoremb}, \cite[Proposition 6.2.4]{bk}).

We recall here that a {\em $T$-cellular variety} is a $T$-variety $X$
equipped with a $T$-stable algebraic cell decomposition. In other
words there is a filtration \be\label{filter} X=Z_1\supseteq
Z_2\supseteq\cdots\supseteq Z_m\supseteq Z_{m+1}=\emptyset\ee where
each $Z_i$ is a closed $T$-stable subvariety of $X$ and
$Z_i\setminus Z_{i+1}=Y_i$ is $T$-equivariantly isomorphic to the
affine space ${\mathbb{C}}^{k_i}$ equipped with a linear action of $T$
for $1\leq i\leq m$. Furthermore, $Y_i$ for $1\leq i\leq m$, are the
Bialynicki-Birula cells associated to a generic one-parameter subgroup
of $T$ (see Definition \ref{cellular} and \cite{u3}). In other words a
$T$-variety is said to be {\em $T$-cellular} if it admits a filtrable
Bialynicki-Birula cellular decomposition with respect to a generic
one-parameter subgroup such that the Bialyncki-Birula cells are
smooth. Any smooth projective complex variety $X$ with $T$-action
having only finitely many $T$-fixed points is $T$-cellular. Indeed for
any projective variety the Bialynicki-Birula cell decomposition is
filtrable (see \cite{BB1}) and since $X$ is smooth the
Bialynicki-Birula cells are smooth (see \cite[Section 3.1]{Br2}).

We let $\mathbb{X}$ be a toroidal $G$-embedding which is $T\times T$-cellular.

Let $\tG_{comp}$ denote a maximal compact subgroup of a factorial
cover $\tG$ of $G$ (see \eqref{exactseq}). Then $\tT_{comp}\subseteq \tG_{comp}$ is the
maximal compact subgroup of the maximal torus $\tT$ of $\tG$.

We consider the action of $\tT$, and hence of $\tT_{comp}$, on $X$
via its projection to $T$. We consider the action of $\tG\times \tG$ (resp. $\tT\times \tT$) and hence of
$\tG_{comp}\times \tG_{comp}$, (resp. $\tT_{comp}\times \tT_{comp}$) on $\mathbb{X}$ (resp. $\bar{G_{ad}}$)
via its projection to $G\times G$ (resp.  $G_{ad}\times G_{ad}$). Under our assumptions $\mathbb{X}$ is 
$\tT\times \tT$-cellular and $X$ is $\tT$-cellular for the induced actions.

In \cite{u1, u2} the $G\times G$ and
$\tG\times \tG$-equivariant algebraic $K$-ring of a regular
embedding of $G$ has been described. There were two types of descriptions. In
\cite[Section 2]{u1}, the description of the $G\times G$-equivariant
$K$-ring of a regular embedding was given in terms of the
Stanley-Reisner system of the associated fan $\mathcal{F}$. In
\cite[Section 3]{u1}, the description of
$G^{ss}\times G^{ss}$-equivariant $K$-ring of the wonderful
compactification $\bar{G_{ad}}$ of $G_{ad}$ was given as an algebra
over the subring generated by the $G^{ss}\times G^{ss}$-equivariant
line bundles on $\bar{G_{ad}}$, where $G^{ss}$ denotes the semisimple
simply connected cover of $G_{ad}$. In \cite{u2}, the second type of description was
extended  to the algebraic $\tG\times \tG$-equivariant
$K$-ring of any regular compactification of $G$, as an algebra over
the $\tT$-equivariant $K$-ring of the toric variety $X^+$ associated
to the fan $\mathcal{F}^+$.  In this description the reason for
considering the action of the factorial cover $\tG$ was firstly to
derive the ordinary algebraic $K$-ring structure from the equivariant
algebraic $K$-ring using the results in \cite{mer} (see \cite[Section
1.2]{u1} and \cite[Section 1.1]{u2}). Another reason to go to the
factorial cover is to use the Steinberg basis for $R(\tT)$ as an
$R(\tG)$-module see \cite[Section 1.2]{u1} and \cite[Section 1.1]{u2}.

%We shall show that the property of being $T\times T$-%cellular for $\mathbb{X}$ is characterized by the %associated toric variety $X=X(\cf)$ being $T$-cellular %(see Section \ref{celltoremb} and Theorem %\ref{equivcellularcriterion}). 

In the present paper our main aim is to extend the description of
\cite{u2} on algebraic $\tG\times \tG$-equivariant $K$-ring of a
regular compactification of $G$, to the topological
$\tG_{comp}\times \tG_{comp}$-equivariant $K$-ring of cellular
toroidal embeddings of $G$. We first give a necessary and sufficient
criterion for a toroidal embedding $\mathbb{X}$ to be
$T\times T$-cellular in terms of the $T$-cellularity of the $T$-toric
variety $X=X(\cf)$ (see Section \ref{celltoremb} and Theorem
\ref{equivcellularcriterion}). Using the above characterization for
cellularity of $\mathbb{X}$ together with the description of the
$T\times T$-fixed points and the $T\times T$-invariant curves of a
toroidal embedding given in Section \ref{invariant curves}, we apply
the GKM type theorem for the equivariant topological $K$-ring of any
$T$-cellular variety proved in Section \ref{GKMsection}, to describe
the $\tT_{comp}\times \tT_{comp}$-equivariant topological $K$-ring of
$\mathbb{X}$. In the setting of topological equivariant $K$-theory,
the results in Section \ref{GKMsection} below replace the results of
Vezzosi and Vistoli \cite{vv} as well as the precise form of
localization theorem in \cite[Theorem 1.3]{u1} on the algebraic
equivariant $K$-ring, which were used earlier for the case of regular
embeddings. We next apply the results in Section \ref{W-action} to
derive that \be\label{W-invariant toroidal} K^0_{\tG_{comp}\times
  \tG_{comp}}(\mathbb{X})\cong K^0_{\tT_{comp}\times
  \tT_{comp}}(\mathbb{X})^{W\times W},\ee since $\mathbb{X}$ is a
$\tG\times \tG$-variety which is $\tT\times \tT$-cellular. Using the
above result together with the description of
$K^0_{\tT_{comp}\times \tT_{comp}}(\mathbb{X})$ in Section
\ref{descriptioneqtop} and Section \ref{$G$-equivariant K-ring}, we
prove our main results namely, the description of the
$\tG_{comp}\times \tG_{comp}$-equivariant $K$-ring of $\mathbb{X}$,
firstly as an algebra over the equivariant topological $K$-ring of
$X^+=X(\cf^+)$ and later as an algebra over the equivariant
topological $K$-ring of $\bar{G_{ad}}$. We further show that the
cellular toroidal embeddings of $G$ are weakly equivariantly formal
for $\tG_{comp}\times \tG_{comp}$-equivariant $K$-theory in Section
\ref{relord} and hence derive the ordinary topological $K$-ring of
$\mathbb{X}$ in Section \ref{ordtoroidal}. Here we remark that in the
setting of topological $K$-theory, the results in Section
\ref{W-action} and Section \ref{relord} replace the results of
Merkurjev \cite{mer}, which were used earlier for describing the
algebraic $\tG\times \tG$-equivariant and ordinary $K$-ring of regular
embeddings.

Here we work with the topological $K$-theory since the
varieties are singular (see Section \ref{algeqkth} below for more
details). In the case of regular embeddings the
$\tG\times \tG$-equivariant algebraic $K$-ring is isomorphic to the
$\tG_{comp}\times \tG_{comp}$-topological $K$-ring (see \cite[Section
1.2.1]{u1}). For the toroidal embeddings of $G$, our results can be
compared with the results on equivariant operational algebraic
$K$-theory in \cite{gon} by using the natural map 
\eqref{eq:transnat} (see Section \ref{algeqkth} below). We remark that our results
also extend the results in \cite{u3} on the description of topological
equivariant $K$-ring of cellular toric varieties to the topological
equivariant $K$-ring of any cellular toroidal embedding.

We fix some further notations before stating our main results.

Let $\Phi$ denote the root system of $(G,T)$.  Let $\Phi^+$ denote the
set of positive roots and $\Delta$ denote its subset of simple roots.  Let $r=\mbox{rank}(G_{ad})=|\Delta|$ is the semisimple rank of $G$. Let $l=\dim(T)=\mbox{rank}(G)$. 

%\{\alpha_1,\ldots, \alpha_r\}

We state below our main results:

\bth\label{Main1} The ring
$K^0_{\tG_{comp}\times \tG_{comp}}(\mathbb{X})$ has the following
direct sum decomposition as a
$K^0_{\tT_{comp}}(X^+)\otimes R(\tG_{comp})$-module: \be\label{Dec2}
K^0_{\tG_{comp}\times \tG_{comp}}(\mathbb{X})=\bigoplus
\prod_{\alpha\in I}(1-e^{\alpha})\cdot K^0_{\tT_{comp}}(X^+)\otimes
R(\tT_{comp})_{I}.\ee This direct sum is a free
$K^0_{\tT_{comp}}(X^+)\otimes R(\tG_{comp})$-module of rank
$|W|$ with basis
$\displaystyle\Big\{\prod_{\alpha\in I}(1-e^{\alpha})\otimes f_{v}:v\in C^{I} ~and
~I\subseteq \Delta\Big\}$

where $C^{I}$ and $f_{v}$ are as defined in Section
\ref{free}. Moreover, the component
$K^0_{\tT_{comp}}(X^+)\otimes 1\subseteq
K^0_{\tT_{comp}}(X^+)\otimes R(\tT_{comp})^{W}$ of the direct sum
decomposition can be identified with the subring of
$K^0_{\tG_{comp}\times \tG_{comp}}(\mathbb{X})$ generated by the
classes of ${\tG_{comp}\times \tG_{comp}}$-equivariant line bundles on
$\mathbb{X}$.

\eeth

\bth\label{Main2}
We have the following isomorphism of
$R(\tG_{comp})\otimes R(\tG_{comp})$-algebras
\[K^0_{\tG_{comp}\times \tG_{comp}}(\mathbb{X})\cong
  K^0_{\tG_{comp}\times
    \tG_{comp}}(\bar{G_{ad}})\otimes_{R(\tT_{comp})}
  K^0_{\tT_{comp}}(X^+)\] where
$K^0_{\tG_{comp}\times \tG_{comp}}(\bar{G_{ad}})$ is an
$R(\tT_{comp})$-algebra via the map which sends $e^{\chi}$ for
$\chi\in X^*(\tT_{comp})$ to $[\mathcal{L}_{\chi}]$ where
$\mathcal{L}_{\chi}$ is the associated
$\tG_{comp}\times \tG_{comp}$-linearized line bundle on $\bar{G_{ad}}$
(see \cite{dp}). Here $K^0_{\tT_{comp}}(X^+)$ is isomorphic as an
$R(\tT_{comp})$-algebra to $PLP(\mathcal{F^+})$ which is the ring of piecewise Laurent polynomial functions on $\cf^+$ (See Section \ref{piecewiselaurent}). \eeth

The following is a brief outline of the paper:

In Section \ref{prelimkth} we recall basic notions on topological
$K$-theory from \cite{segal} and we refer to the link between topological $K$-theory,  algebraic $K$-theory and  operational $K$-theory of \cite{ap}.

In Section \ref{cellular varieties} we recall in detail the 
definition of $T$-cellular varieties and fix some notations.

In Section \ref{eqtopctv} we shall recall the main results from
\cite{u3} on topological equivariant $K$-ring of cellular
$T$-varieties.

In Section \ref{GKMsection} we give a GKM type description of the
$T_{comp}$-equivariant topological $K$-ring of cellular $T$-varieties
(see Theorem \ref{main1}). This is an extension of the result in
\cite{u3} on the $T_{comp}$-equivariant topological $K$-ring of
cellular toric varieties to all cellular $T$-varieties satisfying
Assumption \ref{as}. This result is used in Section
\ref{eqtopcelltoroidal} to describe the $\tT_{comp}\times \tT_{comp}$-equivariant 
$K$-ring of cellular toroidal embeddings.

In Section \ref{W-action} we prove some additional results on the
$T_{comp}$-equivariant and $G_{comp}$-equivariant topological $K$-ring
of a $T$-cellular $G$-variety $X$. In particular, we show the
existence of dot-action of the Weyl group $W$ on $K^0_{T_{comp}}(X)$
and show that $K^0_{T_{comp}}(X)^{W}=K_{G_{comp}}^0(X)$ (see
Proposition \ref{W-action general}). In the case when $X=G/B$ the full
flag variety we compare this with another natural action of $W$ which
is the star action (see \cite{ml}, \cite{kk}). These are topological
$K$-theoretic analogues of the corresponding notions for equivariant
cohomology due to Tymoczko \cite{tym} and equivariant Chow ring due to
Brion \cite{Br2}. In Proposition \ref{W-action general}, we also give
a precise description of the image of $K_{G_{comp}}^0(X)$ as a subring
of $R(T_{comp})^{|X^T|}$. Furthermore, in Section \ref{free} assuming
that $\pi_1(G_{comp})$ is free, we prove that $K^0_{G_{comp}}(X)$ is
free as $R(G_{comp})$-module. We further recall the basis defined by
Steinberg for $R(T_{comp})$ as $R(G_{comp})$-module. These results are
required in Section \ref{eqtopcelltoroidal}, where we study the
$\tG_{comp}\times \tG_{comp}$-equivariant topological $K$-ring of
cellular toroidal embeddings.

In Section \ref{cellulartv} we recall the combinatorial criterion for
a toric variety to be cellular in terms of the associated fan and
recall the description of the $T_{comp}$-equivariant topological
$K$-ring of a cellular toric variety from \cite[Theorem 5.2]{u3} (see
Theorem \ref{tveqktdes1}). We derive Theorem \ref{tveqktdes1} from
Theorem \ref{main1}, since a $T$-cellular toric variety satisfies
Assumption \ref{as}. Furthermore, we also recall from \cite[Theorem 5.6]{u3} the
description of the $T_{comp}$-equivariant topological $K$-ring of a
toric variety as an $R(T_{comp})$-algebra of piecewise Laurent polynomial
functions on the associated fan $\Delta$, denoted $PLP(\Delta)$ (see Section
\ref{piecewiselaurent}).

In Section \ref{celltoremb} we shall give the definition of cellular
toroidal embeddings. Furthermore, we show that that a toroidal
embedding $\mathbb{X}$ is cellular if and only if the associated torus
embedding $X$ of $T$ is cellular (see Theorem \ref{equivcellularcriterion}). More precisely, in Proposition \ref{equivalentcellularcondition} (respectively in Proposition
\ref{equivalentfiltrationcondition}) we show that the Bialynicki-Birula cells of $\mathbb{X}$  are all smooth (resp. the Bialynicki-Birula cell decomposition of $\mathbb{X}$ is filtrable) for a choice of a one-parameter subgroup of $T\times T$, if and only if there exists a one-parameter subgroup of $T$ for which the Bialynicki-Birula cells of $X$ are all smooth (resp. the Bialynicki-Birula cell decomposition of $X$ is filtrable). We note here that the choice of the one-parameter subgroup of $T\times T$ depends on the choice of the one-parameter subgroup of $T$ and vice versa.

In Section \ref{invariant curves} we shall give a precise description
of the $T\times T$-stable curves in $\mathbb{X}$ (see Proposition
\ref{descriptioninvariantcurves}). This extends the description of the $T\times T$-stable curves for a regular embedding of $G$ due to Brion (see \cite{Br1}). We also show in Proposition \ref{toroidalassumption} that a cellular toroidal embedding $X$
satisfies Assumption \ref{as}. This enables us to apply Theorem
\ref{main1} to describe the $\tT_{comp}\times \tT_{comp}$-equivariant
$K$-ring of $\mathbb{X}$ in Section \ref{eqtopcelltoroidal}.

In Section \ref{eqkttvposweyl}, by applying Theorem \ref{tveqktdes1}
we give a description of the $\tT_{comp}$-equivariant topological
$K$-ring of the associated cellular torus embedding $X=X(\cf)$ where
$\cf$ is its associated fan in $N=X_*(T)$ (see Theorem
\ref{main3}). Furthermore, from Theorem \ref{main3} and the $W$-action
on $\cf$(see \eqref{permutation}), we derive the GKM type description of
the $\tT_{comp}$-equivariant topological $K$-ring of the cellular
$T$-toric variety $X^+$ associated to the fan $\cf^+$, which is the
intersection of $\cf$ with the positive Weyl chamber $\mathcal{C}^+$
(see Theorem \ref{cs+}). We also derive the isomorphism
$K_{T_{comp}}^0(X^+)\cong PLP(\cf^+)$ (see Theorem \ref{main2}). Again
these results are used in Section \ref{descriptioneqtop} in the
description of the $\tG_{comp}\times \tG_{comp}$-equivariant $K$-ring
of a cellular toroidal embedding.

The first main aim of Section \ref{eqtopcelltoroidal} is to prove a
GKM type theorem describing the
$\tT_{comp}\times \tT_{comp}$-equivariant $K$-ring of a cellular
toroidal embedding $\mathbb{X}$ as an
$R(\tT_{comp})\otimes R(\tT_{comp})$-subalgebra of
\[K^0_{\tT_{comp}\times \tT_{comp}}(\mathbb{X}^{T\times
    T})\cong\Big(R(\tT_{comp})\otimes R(\tT_{comp})\Big)^{(|W|\times
    |W|)\cdot |\mathcal{F}^+(l)|}.\] More precisely, using the
description of the $T\times T$-fixed points and $T\times T$-stable
curves given in Section \ref{invariant curves}, we give a GKM type
description of $K^0_{\tT_{comp}\times \tT_{comp}}(\mathbb{X})$ (see
Theorem \ref{torus-equivariant K-ring}), by applying Theorem
\ref{main1} on the GKM description of the topological equivariant
$K$-ring of any $T$-cellular variety. This is an extension of the
results on the equivariant topological $K$-ring of a flag variety in
\cite{ml} and that for the equivariant topological $K$-ring of a
cellular toric variety in \cite{u3}.

% describe the $\tT_{comp}\times \tT_{comp}$-equivariant
%topological $K$-ring of a cellular toroidal embedding
%$\mathbb{X}$. More precisely, we
%Firstly as $K^0_{\tT_{comp}}(X^+)\otimes R(\tG_{comp})$-algebra (see
%Theorem \ref{decomp}, Theorem \ref{multstr}) and secondly as an
%algebra over $K^0_{\tG_{comp}\times \tG_{comp}}(\overline{G_{ad}})$
%(see Theorem \ref{relwond} and Corollary \ref{relwondcor}). 

In Section \ref{descriptioneqtop} we prove our main results, namely
the description of $\tG_{comp}\times \tG_{comp}$-equivariant $K$-ring
of a cellular toroidal embedding $\mathbb{X}$. Our main tool in this
section is the GKM description of
$K^0_{\tT_{comp}\times \tT_{comp}}(\mathbb{X})$ along with Proposition
\ref{W-action general}. Using these we first give a description of
$K^0_{\tG_{comp}\times \tG_{comp}}(\mathbb{X})\cong
K^0_{\tT_{comp}\times \tT_{comp}}(\mathbb{X})^{W\times W}$ as a
subring of $\Big(R(\tT_{comp})\otimes R(\tT_{comp})\Big)^{|\cf^+(l)|}$
(see Corollary \ref{$G$-equivariant K-ring} and Corollary
\ref{$G$-equivariant K-ring change of variables}). Following which we
show in Proposition \ref{chain of inclusions} that
$K^0_{\tG_{comp}\times \tG_{comp}}(\mathbb{X})$ has the structure of a
$K^0_{\tT_{comp}}(X^+)\otimes R(\tG_{comp})$-subalgebra of
$K^0_{\tT_{comp}}(X^+)\otimes R(\tT_{comp})$. We then identify
$K^0_{\tT_{comp}}(X^+)\otimes 1$ as the subring generated by
$\tG_{comp}\times \tG_{comp}$-equivariant line bundles on
$\mathbb{X}$. Now, using the Steinberg basis for $R(\tT_{comp})$ as an
$R(\tG_{comp})$-module (see Section \ref{free}) and the description of
$K^0_{\tG_{comp}\times \tG_{comp}}(\mathbb{X})$ in Corollary
\ref{$G$-equivariant K-ring change of variables} and Proposition
\ref{chain of inclusions}, we give a direct sum decomposition of
$K^0_{\tG_{comp}\times \tG_{comp}}(\mathbb{X})$ as a free
$K^0_{\tT_{comp}}(X^+)\otimes R(\tG_{comp})$-module of rank $|W|$ and
give a basis for this free module in Theorem \ref{decomp}. We further
describe the multiplicative structure constants of the basis elements
in Theorem \ref{multstr}. Furthermore, by comparing with the
corresponding direct sum decomposition of
$K^0_{\tG_{comp}\times \tG_{comp}}(\overline{G_{ad}})$ as an
$R(\tT_{comp})\otimes R(\tG_{comp})$-submodule of
$R(\tT_{comp})\otimes R(\tT_{comp})$ (which was earlier proved for the
algebraic equivariant $K$-ring of $\bar{G_{ad}}$ in \cite[Section
3]{u1}), we get the description of
$K^0_{\tG_{comp}\times \tG_{comp}}(\mathbb{X})$ as an algebra over
$K^0_{\tG_{comp}\times \tG_{comp}}(\overline{G_{ad}})$ (see Theorem
\ref{relwond}). Moreover, using the $R(\tT_{comp})$-algebra
isomorphism $K^0_{\tT_{comp}}(X^+)\cong PLP(\cf^+)$ we get Corollary
\ref{relwondcor}. Indeed this
$K^0_{\tG_{comp}\times \tG_{comp}}(\overline{G_{ad}})$-algebra
structure is induced from the canonical
$\tG_{comp}\times \tG_{comp}$-equivariant morphism
$p:\mathbb{X}\ra \bar{G_{ad}}$. These in particular extend the results
in \cite{u2} on $\tG\times \tG$-equivariant algebraic $K$-ring of
regular embeddings of $G$ to
${\tG_{comp}\times \tG_{comp}}$-equivariant topological $K$-ring of a
cellular toroidal embedding of $G$.  Finally in Theorem
\ref{ordkthtoroidal}, we describe the ordinary $K$-ring of a cellular
toroidal embedding as an algebra over the ordinary $K$-ring of
$\bar{G_{ad}}$. Since the description of $K^0(\bar{G_{ad}})$ is known
by the results in \cite[Section 3.1]{u1} we therefore get a complete
description of the topological $K$-ring of a cellular toroidal
embedding.

%Finally in Theorem \ref{ordkthtoroidal},
%we describe the ordinary $K$-ring of a cellular toroidal embedding as
%an algebra over the ordinary $K$-ring of $\bar{G_{ad}}$. Since the
%description of $K^0(\bar{G_{ad}})$ is known by the results in
%\cite{u1} we therefore get a complete description of the topological
%$K$-ring of a cellular toroidal embedding.

\section{Preliminaries on $K$-theory}\label{prelimkth}

Let $X$ be a compact $G_{comp}$-space for a compact Lie group
$G_{comp}$. By $K^0_{G_{comp}}(X)$ we mean the Grothendieck ring of
$G_{comp}$-equivariant topological vector bundles on $X$ with the
abelian group structure given by the direct sum and the multiplication
given by the tensor product of equivariant vector bundles. In
particular, $K^0_{G_{comp}}(pt)$, where $G_{comp}$ acts trivially on
$pt$, is the Grothendieck ring $R(G_{comp})$ of complex
representations of $G_{comp}$. The ring $K^0_{G_{comp}}(X)$ has the
structure of $R(G_{comp})$-algebra via the map
$R(G_{comp})\ra K_{G_{comp}}^0(X)$ which takes
$[V]\mapsto \mathbf{V}$, where $\mathbf{V}=X\times V$ is the trivial
$G_{comp}$-equivariant vector bundle on $X$ corresponding to the
$G_{comp}$-representation $V$. Let $pt$ be a $G_{comp}$-fixed point of
$X$ then the reduced equivariant $K$-ring
${\widetilde K}_{G_{comp}}^0(X)$ is the kernel of the map
$K_{G_{comp}}^0(X)\ra K_{G_{comp}}^0(pt)$, induced by the restriction
of $G_{comp}$-equivariant vector bundles. For $n\in \mathbb{N}$, we
define
$\widetilde{K}^{-n}_{G_{comp}}(X):=\widetilde{K}^0_{G_{comp}}(S^{n}X)$
where $S^{n}X$ is the $n$-fold reduced suspension of $X$. 

If $X$ is locally compact space but not compact we shall consider
$G_{comp}$-equivariant $K$-theory with compact support denoted
$K_{G_{comp},c}^0(X)$. This can be identified with
${\widetilde K}_{G_{comp}}^0(X^+)$ where $X^+$ denotes the one point
compactification of $X$, which is a compact $G_{comp}$-space, where
the point at infinity, which is the base point of $X^+$ is
$G_{comp}$-fixed. We define
${K}_{G_{comp},c}^{-n}(X):=\widetilde{K}^{-n}_{G_{comp}}(X^+)$.

%Since
%$S^n(X^+)\cong (X\times \mathbb{R}^n)^+$, we have
%\[K_{G_{comp},c}^{-n}(X)=K_{G_{comp},c}^{0}(X\times \mathbb{R}^n)\]
% where $G_{comp}$-acts trivially on $\mathbb{R}^n$.

When $X$ is already compact, we define $X^+=X\sqcup {pt}$ as the disjoint
union of $X$ and a base point. In this case we see that
$K^{0}_{G_{comp},c}(X)=\widetilde{K}^0_{G_{comp}}(X^+)=K^{0}_{G_{comp}}(X)$
since
$K^0_{G_{comp}}(X^+)=K_{G_{comp}}^0(X)\oplus K_{G_{comp}}^0(pt)$. In
\cite[Definition 2.8]{segal}, $K^{0}_{G_{comp},c}$ is denoted by
$K^{0}_{G_{comp}}$ without the subscript $c$. We remark here that
$K_{G_{comp},c}^0$ is not a homotopy invariant unlike $K_{G_{comp}}^0$
(see \cite[Proposition 2.3]{segal}).  For example
$K^{0}_{G_{comp}, c}(\mathbb{R}^1)={\widetilde
  K}_{G_{comp}}^0(S^1)=K_{G_{comp}}^{-1}(pt)=0$ whereas
$K^0_{G_{comp},c}(pt)=K^0_{G_{comp}}(pt)=R(G_{comp})$.

We have the equivariant Bott periodicity
$K_{G_{comp}}^{-n}(X)\simeq K_{G_{comp}}^{-n-2}(X)$ given via
multiplication by the Bott element in $K_{G_{comp}}^{-2}(pt)$. (See
\cite{segal}, \cite{at}.) This enables us to define
$K^{n}_{G_{comp}}(X)$ for a positive $n\in \mathbb{Z}$ as
$K^{n-2q}(X)$ for $q\geq n/2$.

For $X$ a locally compact $G_{comp}$-space, $A$ a closed
$G_{comp}$-subspace of $X$ and for $n\in \mathbb{N}$ we define
${K}_{G_{comp}}^{-n}(X, A):=\widetilde{K}_{G_{comp}}^{-n}(X^+, A^+)=\widetilde{K}^0(S^n(X^+/A^+))$ (see
\cite[Definition 2.7, Definition 2.8]{segal}, \cite[Definition 2.2.2]{at}) where $X^+/A+$
is the space obtained by collapsing $A^+$ to a point in $X^+$. 

For $X$ is locally compact $G_{comp}$-space and $A$ is a closed $G_{comp}$-subspace of $X$, the natural inclusion of 
pairs $(X\setminus A, A\setminus A)\hra (X,A)$ induces an isomorphism 
\begin{align*}K_{G_{comp},c}^{-n}(X\setminus A)&\cong \widetilde{K}_{G_{comp}}^0(S^n(X\setminus A)^+)\\ &\cong \widetilde{K}_{G_{comp}}^0(S^{n}(X^+/A^+))\\ &\cong K_{G_{comp}}^{-n}(X,A)\end{align*} because of the identification $(X\setminus A)^+\cong (X^+\setminus  A^+)^+\cong X^+/A^+$ (see \cite[Prop.2.9]{segal}). In particular

Note that $K^{-n}_{G_{comp},c}(\mathbb{C}^q)=\widetilde{K}^{0}_{G_{comp}}(S^n((\mathbb{C}^q)^+))=\widetilde{K}^{0}_{G_{comp}}(S^{n+2q})$ since $(\mathbb{C}^q)^{+}=S^{2q}$. By the equivariant Bott periodicity it follows that when $n$ is even  $K^{-n}_{G_{comp},c}(\mathbb{C}^q)=R(G_{comp})$ and when $n$ is odd $K^{-n}_{G_{comp},c}(\mathbb{C}^q)=0$.

For $A$ a closed subspace of a locally compact space $X$ there is a long exact sequence of $G_{comp}$-equivariant $K$-groups infinite in both
directions given as follows:

{\small \begin{align}\label{l.e.s} \cdots\ra
  {K}_{G_{comp},c}^{-n}(X, A)\ra
  {K}_{G_{comp},c}^{-n}(X) &\ra
  K_{G_{comp},c}^{-n}(A)\ra \nonumber\\ & \ra
  K_{G_{comp},c}^{-n+1}(X, A)\ra
  \cdots \end{align}}

Let $E$ be a $G_{comp}$-equivariant vector bundle on a locally compact space $X$ and $\zeta: X\rightarrow E$ is the zero section. Then \[\uplambda_{-1}(E):=\sum_{k}(-1)^k [\Lambda^k(E)]\in K^0_{G_{comp},c}(X)\] where $\Lambda^k(E)$ denotes the $k$th exterior power of $E$ (see \cite{at}), so that $\Lambda^k(E)=0$ when $k> \dim(E)$.

By the equivariant Thom isomorphism theorem \cite[Proposition 3.2]{segal} there is a class $\uplambda_{E}\in K^0_{G_{comp},c}(E)$ called the {\it Thom class} such that $K_{G_{comp},c}^0(E)$ is a free $K_{G_{comp},c}^0(X)$-module of rank $1$ generated by
$\uplambda_{E}$. Moreover, $\zeta^*(\uplambda_{E})=\uplambda_{-1}(E)$. 
When $E$ is isomorphic to a direct sum of $G_{comp}$-equivariant line bundles $L_1,\ldots, L_n$ then it can be seen that $\displaystyle \uplambda_{-1}(E)=\prod_{i=1}^n (1-[L_i])\in K_{G_{comp}}^0(X)$. 

In particular, when $X=pt$ and $G_{comp}=T_{comp}$ is a compact torus and $E$ is a finite-dimensional $T_{comp}$-representation such that 
$\displaystyle E=\bigoplus_{i=1}^n \mathbb{C}_{\chi_i}$ is a direct sum of $1$-dimensional sub-representations corresponding to the characters $\chi_1,\ldots, \chi_n$ then $\displaystyle \uplambda_{-1}(E)=\prod_{i=1}^n (1-e^{\chi_i})\in K_{T_{comp}}^0(pt)=R(T_{comp})$ and is called the $K_{T_{comp}}$-theoretic equivariant Euler class of $E$ (see \cite[p. 134]{TtD}). Furthermore, the zero section $\zeta$ in this case is the inclusion map of $pt$ in $E$ so that $\zeta^*$ is nothing but the restriction map $K_{T_{comp}}^0(E)\rightarrow K_{T_{comp}}^0(pt)$. Thus the restriction of $\uplambda_{E}$ is $\prod_{i=1}^n (1-e^{\chi_i})\in R(T_{comp})$

\subsection{Algebraic Equivariant $K$-theory}\label{algeqkth}

Let $X$ be an algebraic variety with the action of an algebraic group
$G$. Then $\mathcal{K}^0_{G}(X)$ denotes the Grothendieck ring of
equivariant algebraic vector bundles on $X$ and $\mathcal{K}_0^{G}(X)$
the Grothendieck group of equivariant coherent sheaves on $X$ (see
\cite{Th}, \cite{mer}). The natural map
$\mathcal{K}^0_{G}(X)\ra \mathcal{K}_0^G(X)$ obtained by sending a
class of a $G$-equivariant vector bundle $\mathcal{V}$ on $X$ to the
dual of its sheaf of local sections is an isomorphism when $X$ is
smooth, but not in general. Moreover, when $G$ is a complex reductive
algebraic group and $G_{comp}$ is a maximal compact subgroup of $G$,
then any complex algebraic $G$-variety $X$ is a $G_{comp}$-space. When
$X$ is a smooth compact complex algebraic $G$-variety, we have a natural map
$\mathcal{K}^G_{0}(X)\ra K^0_{G_{comp}}(X)$, obtained by first
identifying $\mathcal{K}^{G}_0(X)$ with $\mathcal{K}_{G}^0(X)$ and
then viewing an algebraic $G$-vector bundle as a topological
$G_{comp}$-vector bundle on $X$ (see \cite[Section 5.5.5]{cg}).

Recall that in \cite[Section 5.5]{cg} we have an alternate notion of
$G$-cellular variety $X$, where $X$ admits a decreasing filtration
\eqref{filter} by $G$-stable closed subvarieties $Z_i$ such that
$Y_i=Z_i\setminus Z_{i+1}$ are complex affine spaces equipped with a
linear $G$-action. However, $Y_i$ are not assumed to be the cells of a
Bialynicki-Birula decomposition of $X$. Thus for $T$-varieties, this
notion of cellular is weaker than our definition (see Definition
\ref{cellular}). It follows from \cite[p. 272]{cg}, that the map
$\mathcal{K}^G_{0}(X)\ra K^0_{G_{comp}}(X)$ is an isomorphism when $X$
is smooth and $G$-cellular.

For singular $G$-varieties there are no natural isomorphisms
$\mathcal{K}^0_{G}(X)\cong \mathcal{K}_0^G(X)$ and
$\mathcal{K}^G_{0}(X)\cong K^0_{G_{comp}}(X)$.

If $T$ is a torus and $X$ a $T$-variety, possibly singular, Anderson and Payne defined (see \cite{ap}) the operational equivariant $K$-ring $opK^0_T(X)$, which is isomorphic to $K^0_T(X)$ when $X$ is smooth. Moreover, it follows from  \cite[Theorem 5.8]{ap} that there exist  natural transformations
\begin{equation}\label{eq:transnat}K^0_{T_{comp}}(X)\to opK^0_T(X)
  \end{equation}
for all  algebraic $T$-varieties $X$.

%they showed that if $X$ is $T$-linear and complete (the cellular varieties are examples of $T$-linear varieties) then there is a natural isomorphism $opK^0_T(X)\simeq Hom_{R(T)}(K^T_0(X),R(T))$ (see \cite[Theorem 6.1]{ap}).

In \cite[Theorem 1.6]{ap} Anderson and Payne prove that the $T$-equivariant  operational $K$-ring of a toric variety $X(\Sigma)$ is isomorphic to the algebra of piecewise exponential functions on $\Sigma$. In \cite[Theorem 5.6]{u3} the analogous result for the topological $T_{comp}$-equivariant $K$-ring of a complete toric variety $X(\Sigma)$ has been proved. 

In \cite[Theorem 5.4]{gon}, Gonzales extended the result \cite[Theorem 1.6]{ap} to any complete $T$-skeletal variety. Our results in Section \ref{GKMsection} are the topological $T_{comp}$-equivariant $K$-theoretic analogue of the results of Gonzales in \cite[Section 5]{gon} for a $T$-cellular variety. In particular, the above natural transformations (\ref{eq:transnat}) are isomorphisms for $T-$cellular varieties. Besides, Gonzales studied the operational equivariant $K$-ring of projective group embeddings in \cite[Section 6]{gon}. In the case of a toroidal group embedding $X$, our description in Theorem \ref{torus-equivariant K-ring} of $\tT_{comp}\times \tT_{comp}$-equivariant topological $K$-ring of $X$ is analogous to the description in \cite[Theorem 6.2]{gon}.

In this article we shall restrict ourselves to the study of topological 
equivariant $K$-theory.

\section{Bialynicki-Birula decomposition and cellular 
  varieties}\label{cellular varieties}

In this Section we recall the notions of a filtrable Bialynicki-Birula
decomposition and a $T$-cellular variety (see \cite[Section 3]{Br2}).

Let $X$ be a normal complex algebraic variety with the action of a
complex algebraic torus $T$. We assume that the set of $T$-fixed
points $X^{T}$ is finite.  A one-parameter subgroup
$\lambda:\mathbb{G}_m\lra T$ is said to be generic if
$X^{\lambda}=X^{T}=\{ x_1,\ldots, x_m\}$.  Let
\be\label{ps}\displaystyle X_{\lambda,x_i}=Y_i=\{ x\in X\mid
\lim_{t\ra 0}\lambda(t)x ~ \mbox{exists and is equal to}~x_i\}.\ee
Then $Y_i$ is locally closed $T$-invariant subvariety
of $X$ and is called the {\it plus} stratum.

\bdefe\label{filtrable} The variety $X$ with the action of $T$ having
finitely many fixed points $\{x_1,\ldots, x_m\}$ is called {\em
  filtrable} if it satisfies the following conditions:

(i) $X$ is the union of its plus strata $Y_i$ for
$1\leq i\leq m$.

(ii) There exists a finite decreasing sequence of $T$-stable closed
subvarieties of $X$ : 
$X=Z_1\supset Z_2 \cdots\supset  Z_{m} \supset Z_{m+1}= \emptyset $
such that $Z_i\setminus Z_{i+1}=Y_i$ for
$1\leq i\leq m$. In particular,
$\displaystyle \overline{Y_i}\subseteq Z_{i}=\bigcup_{j\geq i} Y_j$.
\edefe

\bdefe\label{cellular} Let $X$ be a normal complex algebraic variety
with an action of a complex algebraic torus $T$ such that
$X^T:=\{x_1,\ldots, x_m\}$ is finite. We say that $X$ is {\em
  $T$-cellular} if $X$ is filtrable for a generic one-parameter
subgroup $\lambda$ and if moreover, each plus stratum $Y_i$ is
$T$-equivariantly isomorphic to a complex affine space
$\mathbb{C}^{n_i}$ on which $T$-acts linearly for $1\leq i\leq
m$. \edefe

%We were also motivated by
%the definition of (possibly singular) varieties with {\it good
%  decompositions} due to Carrell and Goresky (see \cite{cago} and
%.  The $T$-cellular varieties
%which we consider are varieties which admit good decompositions in the
%above sense.

\brem\label{carrell}  Carrell and Goresky (see \cite{cago} and \cite[Section 4.3, Theorem 4.13]{car}) study possibly singular varieties having {\it good decompositions}, that is filtrable Bialynicki-Birula cellular decomposition and compute the integral homology for such varieties. The
$T$-cellular varieties that we consider are varieties that admit good
decompositions in the above sense.\erem

\section{Equivariant topological $K$-ring of cellular
          varieties}\label{eqtopctv}

Let $T\simeq (\mathbb{C}^*)^n$ and
$T_{comp}\simeq (S^1)^n\subseteq T$.

We shall assume that the $T$-cellular variety $X$ is compact, which in turn implies that $Z_i$ is compact for $1\leq i\leq m$ and $Z_m=Y_m=\{x_m\}$.

\bth\label{k-thcellular} Let $X$ be a $T$-cellular variety (see
Definition \ref{cellular}). Then $K_{T_{comp}}^0(X)$ is a free
$R(T_{comp})$-module of rank $m$-which is the number of
cells. Furthermore, we have $K_{T_{comp}}^{-1}(X)=0$.
\eeth

\begin{proof}
  By  \cite[Section 2]{segal}, \cite{at}) (also see Section \ref{prelimkth}) we have the following long exact sequence
  of $T_{comp}$-equivariant $K$-groups which is
  infinite in both directions: \begin{align}\label{les} \cdots\ra
    {K}^{q}_{T_{comp}}(Z_i,Z_{i+1})\ra {K}_{T_{comp}}^{q}(Z_i)&\ra
                                                                K^{q}_{T_{comp}}(Z_{i+1})\ra \nonumber\\
                                                              &\ra
                                                                K^{q+1}_{T_{comp}}(Z_i,Z_{i+1})\ra
                                                                \cdots \end{align}
                                                                for
                                                                $1\leq
                                                                i\leq
                                                                m$ and
                                                                $q\in
                                                                \mathbb{Z}$.

Moreover, by \cite[Proposition 2.9]{segal}(also see Section \ref{prelimkth}) we have
\begin{align*}{K}^q_{T_{comp}}(Z_i,Z_{i+1}) &\cong K^q_{T_{comp},c }(Y_i)\\
                                         &\cong
                                           K_{T_{comp}}^{q-2n_i}(x_i) \end{align*}
                                           where
                                           $Z_i\setminus Z_{i+1}=Y_i=\mathbb{C}^{n_i}$
                                           for $1\leq i\leq m$ . Thus
                                           when $q$ is even
                                           $K^{q}_{T_{comp}}(Z_i,Z_{i+1})=K^0_{T_{comp}}(pt)=R(T_{comp})$,
                                           and when $q$ is odd
                                           $K^{q}_{T_{comp}}(Z_i,Z_{i+1})=0$. The last isomorphism above is also the equivariant Thom isomorphism for the $T_{comp}$-equivariant vector bundle $Y_i$ over the point $x_i$ (see \cite[Proposition 3.3]{segal})

                                           Moreover, since $Z_m=Y_m=\{x_m\}$
                                           and $Z_{m+1}=\emptyset$ we
                                           have
                                           $K^0_{T_{comp}}(Z_m)=
                                           K^0_{T_{comp}}(pt)=R(T_{comp})$
                                           and
                                           $K^{-1}_{T_{comp}}(Z_m)=K^{-1}_{T_{comp}}(pt)=0$.

    Now, by descending induction on $i$, we suppose that
    $K_{T_{comp}}^0(Z_{i+1})$ is a free $K^0_{T_{comp}}(pt)$-module
    for $1\leq i\leq m$ of rank $(m-i)$ and
    $K_{T_{comp}}^{-1}(Z_{i+1})=0$.  We can start the induction since
    $K^0_{T_{comp}}(Z_m)=K_{T_{comp}}^0(pt)$ and
    $K^{-1}_{T_{comp}}(Z_m)=0$.

    Thus from \eqref{les} we get the following split short exact
    sequence of $K_{T_{comp}}^0(pt)$-modules \be \label{eq1} 0\ra
    K_{T_{comp}}^0(Z_i, Z_{i+1})\ra K_{T_{comp}}^0(Z_i) {\ra} K_{T_{comp}}^0( Z_{i+1})\ra 0 \ee for
    $1\leq i\leq m$.

    From \eqref{eq1} and induction, $K_{T_{comp}}^0(Z_i)$ is a free
    $R(T_{comp})$-module of rank $m-i+1$. In particular, since $X=Z_1$, it follows that
    $K_{T_{comp}}^0(X)$ is a
    free $R(T_{comp})$-module of rank $m$.

    Now, since $K^{-1}_{T_{comp}}(Z_i,Z_{i+1})=0$ as shown above, and 
    $K^{-1}_{T_{comp}}(Z_{i+1})=0$ by 
    induction assumption,
    it follows from \eqref{les}, that
    $K^{-1}_{T_{comp}}(Z_i)=0$. In particular, since $X=Z_1$, 
    $K_{T_{comp}}^{-1}(X)=0$.

\end{proof}

Recall that $X^{T}=\{x_1,\ldots, x_m\}$ and $Z_i^{T}=\{x_i,\ldots, x_m\}$ for each $1\leq i\leq m$. Let $\iota:X^{T_{comp}}=X^{T}\hra X$ denote the inclusion of the set of
$T$-fixed points in $X$. Let $\iota_i:Z_i^{T_{comp}}=Z_i^{T}\hra Z_i$ denote the inclusion of the set of $T$-fixed points for $1\leq i\leq m$.

For a $T$-cellular variety $X$, we prove the following proposition.

\bpropo\label{localization} The canonical restriction
map
\[ K^0_{T_{comp}}(X)\stackrel{\iota^*}{\lra} K^0_{T_{comp}}(X^{T_{comp}})\cong R(T_{comp})^m\] is an
injective map of $R(T_{comp})$-algebras.  \epropo
\begin{proof} Consider the following commuting diagram where the first row is
   \eqref{eq1}, and the vertical arrows are the maps
   $ \iota^*_{i}$ and $\iota^*_{i+1}$ induced by
   restriction to the $T$-fixed points of $Z_i$ and $Z_{i+1}$
   respectively:

\begin{center}
	\begin{tikzcd}[cramped,column sep=tiny]
          0\ar[r] & R(T_{comp})=K_{T_{comp}}^0(x_i)\arrow[r]
          \arrow[d,-,double equal sign distance,double]
          &K_{T_{comp}}^0(Z_i)\arrow[r] \arrow[d]
          &K_{T_{comp}}^0(Z_{i+1}) \arrow[d]\ar[r] &0\\
          0\ar[r] & R(T_{comp})=K_{T_{comp}}^0(x_i)\arrow[r] &
          K^0_{T_{comp}}(Z_i^{T}) \arrow[r] &
          K^0_{T_{comp}}(Z_{i+1}^{T})
          \ar[r] &0 \\
	\end{tikzcd}
\end{center}

The proof of injectivity of
\[\iota^*_{i}: K_{T_{comp}}^0(Z_i)\ra
K^0_{T_{comp}}(Z_i^{T})=
R(T_{comp})^{m-i+1}\] follows by descending
induction on $i$ and diagram chasing.  We can start the induction
since $Z_m=\{x_m\}$ and
$\iota^*_{m}=id_{R(T_{comp})}$. Since $Z_1=X$
the proposition follows.
  
  \end{proof}

  \brem\label{alternate} The above proposition can alternately be
  proved using the localization theorem \cite[Prop.4.1]{segal} and
  Theorem \ref{k-thcellular} (see \cite[Cor. 4.3]{u3}).  \erem

\subsection{\boldmath GKM theory of $X$}\label{GKMsection}

In this section we give a GKM type description for $K^0_{T_{comp}}(X)$, as a
$K^0_{T_{comp}}(pt)=R(T_{comp})$-subalgebra of $R(T_{comp})^{m}$, for a complete $T$-cellular variety $X$.

Recall that (see \cite{cg}) we have the isomorphism \be\label{iso3'}
R(T)\simeq R(T_{comp}). \ee 

For a $T$-cellular variety $X$, with $\mbox{dim}(X) =n$, we have a decreasing sequence of $T$-stable closed
subvarieties  such that
$Z_i\setminus Z_{i+1}=Y_i\simeq \bc^{n-k_i}$. Moreover,
$x_i\in Y_i\subseteq Z_i$, $1\leq i\leq m$. Since
$Y_i\simeq\mathbb{C}^{n-k_i}$ are $T$-stable we have a
$T$-representation $\rho_i=(V_i,\pi_i,x_i)$, which can alternately be
viewed as a $T$-equivariant complex vector bundle over the $T$-fixed point $x_i$, for $1\leq i\leq m$.  

Furthermore, there exist some $\chi_i^l$, characters of $T$, such that 

\be \label{subrep}
V_i\cong \bigoplus_{l=1}^{n-k_i} \mathbb{C}_{\chi_{i}^l} \ee

 as a
$T$-representation, where $\mathbb{C}_{\chi_{i}^l}$ is the
$1$-dimensional $T$-sub-representation associated to the character
$\chi_{i}^l$.

\bdefe\label{outin}
We shall call a curve $C$ incident at a $T$-fixed point $x_i$ as
outgoing at $x_i$ if $C$ joins $x_i$ and $x_j$ with $j>i$,
otherwise we call $C$ as incoming at $x_i$.\edefe

\brem For  $C$  a $T-$stable irreducible curve, outgoing  at $x_i$ means that the plus strata $C\cap Y_i$ is open in $C$ and incoming at $x_j$ means that the plus stata $C\cap Y_j = x_j$ is just a point.

\erem
\bnot\label{outnot}
For each $x_i\in X^T$, we denote by
\[\mathcal{C}^0_i\] the set of outgoing
$T$-stable curves in $X$ at $x_i$. For each $T$-stable irreducible curve $C$ in $X$, we denote by $\chi_C$ the character by which  $T$ acts on $C$.\enot

We make the following assumptions.

\bas\label{as} 

\begin{enumerate} 

\item In $X$, every $T$-stable curve joining a pair of
  $T$-fixed points is isomorphic to $\mathbb{P}^1_{\mathbb{C}}$.

\item For each $x_i\in X^T$, there are precisely $|\mathcal{C}^0_i|=n-k_i$ outgoing curves at $x_i$ which are
  in one to one correspondence with the one-dimensional
  $T$-subrepresentations of $V_i$.

\item For each $x_i\in X^T$, the characters $\chi_{C} $, for $C\in \mathcal{C}^0_i$, are pairwise linearly independent in the lattice
$M=Hom(T,\mathbb{C}^*)$.  \end{enumerate}\eas

Let $\mathcal{A}$ be the set of all
$(a_i)_{1\leq i\leq m}\in R(T_{comp})^m$ such that
$a_i\equiv a_j \pmod{(1-e^{\chi_{i,j}})}$ whenever the $T$-fixed
points $x_i \neq x_j$ lie in a closed $T$-stable irreducible curve
$C_{i,j}\cong \mathbb{P}^1_{\mathbb{C}}$ in $X$ and $T$ (and hence
$T_{comp}$) acts on $C_{i,j}$ through the character
$\chi_{i,j}$.\footnote{For $i<j$, there could be several different $T$-stable irreducible curves which pass by $x_i$ and $x_j$ but in the cases we are interested in (cf Remark \ref{rem:onecurve}), there will be at most one irreducible $T-$stable curve which contain $x_i$ and $x_j$. So there will be no ambiguity for labelling curves by $C_{i,j}$ and the corresponding characters by $\chi_{i,j}$.} Moreover, $\mathcal{A}$ is a $R(T_{comp})$-algebra where
$R(T_{comp})$ is identified with the subalgebra of $\mathcal{A}$
consisting of the diagonal elements $(a,a,\ldots, a)$.

Let $\iota:X^{T}\hra X$ be the inclusion of the set of $T$-fixed
points in $X$.

The following lemma is a refinement of Proposition \ref{localization}
for a $T$-cellular toric variety $X$. For a similar result
when $X$ is a full flag variety see \cite[Lemma 2.3]{ml}.

\blem\label{localgkm} Let $X$ be a $T$-cellular variety. Then we have
the following inclusion of $R(T_{comp})$-subalgebras
$\iota^*(K^0_{T_{comp}}(X))\subseteq \mathcal{A}$ of $R(T_{comp})^m$.  \elem
\begin{proof} Let the $T$-fixed points $x_i$ and $x_j$ lie in the
  closed $T$-stable irreducible curve $C_{i,j}$ and $T$ (and hence
  $T_{comp}$) act on $C_{i,j}$ through the character $\chi_{i,j}$. The
  composition
  \[K^0_{T_{comp}}(X)\stackrel{\iota^*}{\ra} K^0_{T_{comp}}(X^{T})\ra K^0_{T_{comp}}(\{x_i,x_j\})\] equals
  \[K^0_{T_{comp}}(X)\ra
    K^0_{T_{comp}}(C_{i,j})\stackrel{\iota_{C_{i,j}}^*}{\ra}
    K^0_{T_{comp}}(\{x_i,x_j\}).\]

  Furthermore, by our assumption $C_{i,j}$ is smooth and isomorphic to
  $\mathbb{P}_{\mathbb{C}}^1$. We now claim that the image of
  $\iota_{C_{i,j}}^*$ consists of pairs of elements
  $(f,g)\in R(T_{comp})\oplus R(T_{comp})$ such that
  $f\equiv g\pmod{(1-e^{\chi_{i,j}})}$.

  The proof of this fact follows by the same arguments as can be found
  in \cite[proof of Lemma 2.3, p. 322]{ml} and \cite[Theorem 1.3,
  p. 274]{u1} so we do not repeat them (also see \cite[Section
  3.4]{Br2}).

  Thus if $x\in K^0_{T_{comp}}(X)$ then $\iota^*(x)_{i}-\iota^*(x)_j$ is
  divisible by $(1-e^{\chi_{i,j}})$. Hence the lemma.

\end{proof}

\bth\label{main1} Let $X$ be a complete $T$-cellular variety
satisfying assumptions \ref{as}.  The ring $K^0_{T_{comp}}(X)$ is isomorphic to
$\mathcal{A}$ as an $R(T_{comp})$-algebra.

\eeth \begin{proof} By Corollary \ref{localization} and Lemma
  \ref{localgkm} we know that the image of
  $\iota^*(K^0_{T_{comp}}(X))\subseteq \mathcal{A}$. It remains to
  show that $\iota^*(K^0_{T_{comp}}(X))\supseteq \mathcal{A}$. We
  shall follow the methods in \cite[Section 2.5]{ml} and \cite[Theorem 5.2]{u3}.

  Recall that for $1\leq i\leq m$ there are outgoing $T$-stable irreducible curves at $x_i$ :
  \[C_i^j,\, 1\le j \le n-k_i.\]
  
  For each $j$, the curve $C_i^j$ joins  $x_i$ and $x_{i_j}$, with $i<i_j$, and $T$ acts on it by a character $\chi_i^j$.

  Since $\chi_{i}^j$ is non-zero for every $1\leq j\leq n-k_i$, the
  $K$-theoretic equivariant Euler class
  $1-e^{\chi_{i}^j}$ is a non-zero divisor in the
  unique factorization domain $R(T_{comp}) =\mathbb{Z}[e^u : u\in M]$.
  Moreover, since $\chi_i^j$ and $\chi_i^{j'}$ are linearly
  independent, it follows that $1-e^{\chi_i^j}$ and
  $1-e^{\chi_i^{j'}}$ are relatively prime in $R(T_{comp})$.

By the Thom isomorphism $K^0_{T_{comp}}(V_i)$ is a free
$R(T_{comp})$-module on one generator $g_i$ which restricts under the
restriction map $K^0_{T_{comp}}(V_i)\lra K^0_{T_{comp}}(x_i)$ to
$\uplambda_{-1}(V_i)=\displaystyle\prod_{j=1}^{n-k_i}
(1-e^{\chi_i^j})$ (see \cite{segal} and Section \ref{prelimkth}).

Recall from \eqref{eq1} that for every $1\leq i\leq m$ we have a split
exact sequence of $R(T_{comp})$-modules
\be\label{splitses} 0\ra
  K^0_{T_{comp}}(Y_i)\stackrel{j_i}{\ra}K^0_{T_{comp}}(Z_i)\stackrel{k_i}{\ra}K^0_{T_{comp}}(Z_{i+1})\ra
  0\ee where $K^0_{T_{comp}}(Y_i)=K^0_{T_{comp}}(Z_i,Z_{i+1})$.  These
induce the following maps induced by inclusions
\[\begin{array}{lllllllll}
  K^0_{T_{comp}}(x_1) & &K^0_{T_{comp}}(x_2) &  & &K^0_{T_{comp}}(x_{m-1}) & &K^0_{T_{comp}}(x_m)\\
  ~~~~\uparrow~~~&  &~~~~\uparrow~~~&  & &~~~~\uparrow~~~& &~~~~\uparrow ~~~\\
  K^0_{T_{comp}}(Z_1=X)& \stackrel{k_1}{\ra}&K^0_{T_{comp}}(Z_2)\stackrel{k_2}{\ra}&\cdots&\ra&K^0_{T_{comp}}(Z_{m-1})&\stackrel{k_{m-1}}{\ra}&K^0_{T_{comp}}(Z_m)\\
  ~~~~\uparrow j_{1}~~~& &~~~~\uparrow j_2~~~& & &~~~~\uparrow j_{m-1}~~~& &~~~~\uparrow j_{m}~~~\\
 K^0_{T_{comp}}(Y_1)&  &K^0_{T_{comp}}(Y_2)&  & & K^0_{T_{comp}}(Y_{m-1})& &K^0_{T_{comp}}(Y_m)  
\end{array}\]

In the above diagram $k_i$ are all surjective and $j_i$ are all injective maps.

Thus we can choose $f_i\in K^0_{T_{comp}}(X)$ such that
$k_{i-1}\circ\cdots k_1(f_i)=j_i(g_i)$. Thus
$\displaystyle \iota^*(f_i)_{i}=\prod_{j=1}^{n-k_i} (1-e^{\chi_{i}^j})$
and $\iota^*(f_i)_{{l}}=0$ if $l>i$ since $k_i\circ j_i=0$.

Let $a\in \mathcal{A}$. We prove by induction on $i$ that if $a_{l}=0$
for $l\geq i$, then $a\in \iota^*(K^0_{T_{comp}}(X))$. This will complete the proof of the theorem.

For $i=1$, since $a=0$ the claim is obvious.

Assume by induction that $a_{l}=0$ for $l\geq i+1$. We know by the
above arguments that there are $T$-stable irreducible curves $C_i^j$ in $X$, $1\le j \le n-k_i$, outgoing at $x_i$, joining $x_i$ to $x_{i_j}$, with $i<i_j$, and on
which $T$ acts through the character $\chi_{i}^j$. Since $i_j>i$ by
induction hypothesis $a_{i_j}=0$ for all $1\leq j\leq n-k_i$. Now,
since $a\in \mathcal{A}$ we have $(1-e^{\chi_{i}^j})$ divides
$a_i-a_{i_j}=a_i$ for $1\leq j\leq n-k_i$. This implies that
$\displaystyle\prod_{1\leq j\leq n-k_i} (1-e^{\chi_{i}^j})$ divides
$a_i$ since $R(T_{comp})$ is a unique factorization domain and the
factors are relatively prime. Let
\[a_i=c_i\cdot \prod_{1\leq j\leq n-k_i}(1-e^{\chi_{i}^j})\]
where $c_i\in R(T_{comp})$. Thus
$a-\iota^*(c_i\cdot f_i)\in \mathcal{A}$ by Lemma \ref{localgkm}. Also
$(a-\iota^*(c_i\cdot f_i))_{i}=0$. Thus we get
$(a-\iota^*(c_i\cdot f_i))_{l}=0$ for $l\geq i$. By induction there
exists $q\in K^0_{T_{comp}}(X)$ such that
$\iota^*(q)=a-\iota^*(c_i\cdot f_i)$. Thus we get
$a=\iota^*(q+c_i\cdot f_i)$. Hence the theorem.
\end{proof}

\subsection{\boldmath The action of $W$ on $K_{T_{comp}}^0(X)$}\label{W-action}

We shall follow the notations in Section \ref{Introduction}.

Let $X$ be a $G$-variety that is $T$-cellular. Suppose that $X$
satisfies assumptions \ref{as}.

Let $w\in W$ and $w=nT$ for $n\in N(T)$. Consider the action of
$n^{-1}$ from the left on the $G$-space $X$. This action is
$T$-equivariant where we give the usual $T$-action on the left and the
$T$-action twisted by $w$ on the right. Indeed, $\forall x \in X,\, \forall t \in T,\, n^{-1}\cdot t\cdot x=t\# (n^{-1} \cdot x)$ where
$\forall t \in T,\, \forall y \in X,\, t\# y=n^{-1}tn\cdot y$. This can be seen to be independent of the
choice of $n\in N(T)$ such that $w=nT$.

The action by $n^{-1}$ induces morphism
$K_{T_{comp}}^0(X)\stackrel{(n^{-1})^{*}}{\ra} K_{T_{comp}}^0(X)$ by
pull back. Suppose $n'=nt$ then $n'^{-1}=t^{-1}\cdot n^{-1}$ and
$(n'^{-1})^{*}=(n^{-1})^*\cdot (t^{-1})^*$. Since $(t^{-1})^*$ induces
identity map on $K_{T_{comp}}^0(X)$, $(n'^{-1})^*=(n^{-1})^*$.

Thus we have a left action of $W$ on $K_{T_{comp}}^0(X)$ defined by
$(w^{-1})^*=(n^{-1})^*$ for $n\in N(T)$ such that $w=nT$. We define, for $\mathcal{V}$, a $T_{comp}$-equivariant vector bundle on $X$, 
\be\label{w-action} w([\mathcal{V}]):=[(w^{-1})^*(\mathcal{V})] \ee 

Note that $(w^{-1})^*(\mathcal{V})\mid_{x}=\mathcal{V}\mid_{n^{-1}\cdot x}$
and the linear map
$(w^{-1})^*(\mathcal{V})\mid_{x}\ra (w^{-1})^*(\mathcal{V})\mid_{t\cdot x}$ induced
by $t\in T_{comp}$ is the linear map
$\mathcal{V}\mid_{n^{-1}\cdot x}\ra \mathcal{V}\mid_{t\#
  n^{-1}\cdot x}$ induced by $t'=n^{-1}\cdot t\cdot n$.

For a
  $T_{comp}$-representation $V$ and $w\in W$ we define the $T_{comp}$-representation $w(V)$ to be the the vector space $V$ with the
  $T_{comp}$-action given by
  $t\# v:=(n^{-1}\cdot t\cdot n)\cdot v$.

  In particular, when $x\in X^{T}$ then it can be seen that
  $n^{-1}\cdot x\in X^{T}$ and $(w^{-1})^*(\mathcal{V})\mid_{x}$ is the
  $T_{comp}$-representation
  $w\Big(\mathcal{V}\mid_{n^{-1}\cdot x}\Big)$. 

  When $\mathcal{V}=X\times V$ is a trivial $T_{comp}$-equivariant
  vector bundle on $X$, then $\mathcal{V}\mid_{x}=\{x\}\times V$ where
  $V$ is a $T_{comp}$-representation space. Thus we can see that
  $(w^{-1})^*(X\times V)\cong X\times w(V)$. In particular, the
  $W$-action on $K_{T_{comp}}^0(X)$ preserves the
  $R(T_{comp})$-algebra structure on $K_{T_{comp}}^0(X)$ obtained by
  sending an isomorphism class $[V]$ of a $T_{comp}$-representation
  $V$ to $[X\times V]$ with the diagonal $T_{comp}$-action. 

  Moreover, when $\mathcal{V}$ is a $G_{comp}$-equivariant vector
  bundle then $(n^{-1})^*(\mathcal{V})\cong \mathcal{V}$ for $n\in N(T)$
  such that $nT=w$. Since this is an isomorphism of
  $G_{comp}$-equivariant vector bundles, it is also an isomorphism of
  $T_{comp}$-equivariant vector bundles via the forgetful map
  $\phi: K_{G_{comp}}^0(X)\ra K_{T_{comp}}^0(X)$. Hence the induced
  action of $w\in W$ on $K_{G_{comp}}^0(X)$ is given by
  $w([\mathcal{V}])=[\mathcal{V}]$. Thus under $\phi$ the image of
  $K_{G_{comp}}^0(X)$ lies in the subring of $W$ invariants $K_{T_{comp}}^0(X)^{W}$.

  Moreover, if $x\in X^{T}$ then $n^{-1}\cdot x=n'^{-1}\cdot x$ whenever $nT=n'T$. Thus
  we have an action of $W$ on $X^{T}$ given by
  $w^{-1}\cdot x:=n^{-1}\cdot x$ where $w=nT$ for $n\in N(T)$. We
  shall denote the orbits of $X^{T}$ under $W$-action by $X^{T}/W$.

%There is a $W$-action on $X^{T}$ defined as follows. Let $x\in
%X^T$. Define $w\cdot x:=n\cdot x$ where $n\in N(T)$ and
%$nT=w$. Suppose $w=n'T$ for $n'\in N(T)$ then $n'=nt$ for some
%$t\in T$. Thus $n'\cdot x=nt\cdot x=n\cdot x$ and the action is well
%defined.  Furthermore, $t\in T$ let $ntn^{-1}=t'\in T$. Note that
%$t\cdot (w \cdot x)=t\cdot n\cdot x=n\cdot t'\cdot x=n\cdot x=w\cdot
%x$.

  Note that $W$ acts on $T$-stable curves in $X$ as follows. We define
  $w^{-1}\cdot C:=n^{-1}\cdot C$ where $n\in N(T)$ is such that
  $w=nT$. Since $C$ is $T$-stable the action is well defined as it is
  independent of the representative for $w$ chosen in
  $N(T)$. Moreover, $n\in G$ so acts as an isomorphism of $X$ so that for each $T$-stable irreducible curve $C$ in $X$,
  $n^{-1}\cdot C$ is also a $T$-stable irreducible curve in $X$.

  In particular, if $C$ joins the fixed points $x$ and $x'$ then
  $w^{-1}\cdot C$ joins the $T$-fixed points $w^{-1}\cdot x$ and
  $w^{-1}\cdot x'$. Furthermore, if $T$ acts on $C$ via the character
  $\chi$ then $T$ acts on $w^{-1}\cdot C$ via the character $w^{-1}(\chi)$.

We have the following proposition.

\bpropo\label{W-action general} There is a canonical $W$-action on $K^0_{T_{comp}}(X^{T})$ which
restricts to the $W$-action on $K^0_{T_{comp}}(X)$. Moreover,
$K^0_{T_{comp}}(X)^{W}\cong K_{G_{comp}}^0(X)$ (see \cite{ml}). Under
$\iota^*$, $K_{G_{comp}}^0(X)$ maps onto the subring $\mathcal{B}$ of
$\mathcal{A}$ where $\mathcal{B}$ is the set of all
$(a_{{x}})_{{x}\in X^{T}}\in R(T_{comp})^{|X^{T}|}$ such that
$a_{x}\equiv a_{x'} \pmod{(1-e^{\chi})}$ whenever the $T$-fixed points
$x$ and $x'$ lie in a closed $T$-stable irreducible curve $C$ in $X$
and $T$ (and hence $T_{comp}$) acts on $C$ through the character
$\chi$ and $a_{w.x}=w.a_x\in R(T_{comp})$ for all $x\in X^T$ and all $w\in W$.  \epropo

\begin{proof} 
  For all $w \in W$, the above defined action of $w^{-1}\in W$ on $X^{T}$ induces the
  following left action
  $(w^{-1})^*:K_{T_{comp}}^0(X^{T}){\ra} K^0_{T_{comp}}(X^{T})$ by pull
  back. By identifying
  $\displaystyle K_{T_{comp}}^0(X^{T})\cong \bigoplus_{x\in X^{T}}
  K_{T_{comp}}^0(x)$ this action is given by : $\forall w\in W,\, w\cdot (a_{x})=(b_x)$
  where $b_{x}=w\cdot a_{w^{-1}\cdot x}$.

  Let $(a_{x})\in K_{T_{comp}}^0(X)$. Thus
  $a_{x}\equiv a_{x'}\pmod {1-e^{-\chi}}$ whenever $x$ and $x'$ are
  joined by a $T$-stable irreducible curve $C$ on which $T$ acts via the character
  $\chi$. Recall that whenever $x$ and $x'$ are joined by a $T$-stable
  curve $C$ on which $T$ acts via the character $\chi$ then
  $w^{-1}\cdot x$ and $w^{-1}\cdot x'$ are joined by a $T$-stable
  curve $w^{-1}\cdot C$ on which $T$ acts via $w^{-1}(\chi)$. Thus
  \be\label{transequiv} a_{w^{-1}\cdot x}\equiv a_{w^{-1}\cdot x'}\pmod
  {1-e^{-w^{-1}(\chi)}}.\ee Since $b_{x}=w\cdot a_{w^{-1}\cdot x}$, \eqref{transequiv} implies that
  $b_{x}\equiv b_{x'} \pmod {1-e^{-\chi}}$. Thus $(b_{x})\in K_{T_{comp}}^0(X)$.

  This induced action on $K_{T_{comp}}^0(X)$ by restricting the above
  $W$-action on $K_{T_{comp}}^0(X^{T})$ via $\iota^*$, is the same
  action given by the pull backs $(w^{-1})^*$ of $T_{comp}$-equivariant
  vector bundles defined in \eqref{w-action}, $w\in W$.

We have the canonical map
\be\label{characteristic}\Psi: R(T_{comp})\otimes_{R(G_{comp})}K_{G_{comp}}^0(X)\ra
  K_{T_{comp}}^0(X)\ee defined by
$\Psi([{V}]\otimes 1):=[X\times {V}]$ where
$T_{comp}$ acts diagonally and
$\Psi(1\otimes [\mathcal{W}])=\phi([\mathcal{W}])$. Since
$K_{T_{comp}}^0(X)$ is a free $R(T_{comp})$-module by Theorem
\ref{k-thcellular}, by \cite{ml} $\Psi$ is an isomorphism of
$R(T_{comp})$-modules. Moreover, $\Psi$ is $W$-equivariant with
respect to the $W$-action on $K_{T_{comp}}^0(X)$ defined above and the
$W$-action on the left hand side is given by
$w\cdot \sum_{i} [V_i]\otimes [\mathcal{W}_i]:=\sum_i [w(V_i)]\otimes
[\mathcal{W}_i]$. Since $R(T_{comp})^W=R(G_{comp})$, by taking
$W$-invariants on both sides of the isomorphism $\Psi$, we get that
the forgetful map $\phi$ maps $K^0_{G_{comp}}(X)$ isomorphically onto
the subring $K_{T_{comp}}^0(X)^{W}$ of $W$-invariants in
$K_{T_{comp}}^0(X)$.

Thus under the restriction map
$\iota^*: K^0_{T_{comp}}(X)\ra K^0_{T_{comp}}(X^{T})$,
$K^0_{G_{comp}}(X)\cong K_{T_{comp}}^0(X)^{W}$ maps to the subring
consisting of $(a_{x})\in \mathcal{A}$ such that
$w\cdot a_{w^{-1}\cdot x}=a_x$ for $w\in W$ and $x\in X^{T}$ which is
the subring
  \[\{ (a_{x})\in \mathcal{A}~\mid ~ a_{w^{-1}\cdot x}=w^{-1}\cdot
    a_{x} ~\mbox{for}~ w\in W~ \mbox{and}~ x\in X^{T}\}.\] Now, if
  $a_x-a_{x'}$ is divisible by $(1-e^{-\chi})$ then
  $w^{-1}\cdot a_x-w^{-1}\cdot a_{x'}$ is divisible by
  $(1-e^{-w^{-1}(\chi)})$. Hence it follows that the above subring of
  $\mathcal{A}$ is isomorphic to $\mathcal{B}$.
\end{proof}

\subsubsection{Comparison between the Star and Dot Action in
  topological equivariant $K$-theory}\label{compstardot}

The above action of $W$ on $K_{T_{comp}}^0(X)$ is the topological
$K$-theoretic analogue of the {\it dot} action in equivariant
cohomology defined by Tymoczko (see \cite[Section 3]{tym}) and earlier for the
equivariant Chow ring by Brion (see \cite[Section 6.2]{Br2}). In the case when
$X=G/B$ there is another action of $W$ on $K_{T_{comp}}^0(X)$ (see
\cite{ml}, \cite{kk}) which is induced by the action of $W$ on $G/B$
from the right given by $gB\cdot w:=gnB$ where $n\in N(T)$ is such
that $nT=w$. The above action of $W$ commutes with the $T$-action from
the left and hence induces a left action of $W$ on $G/B$ via pull back
of $T_{comp}$-equivariant vector bundles via the map
$gB\mapsto gB\cdot w$.

For the above action of $W$ on $G/B$, the restriction map $\iota^*$ to
the $T$-fixed points of $G/B$, is $W$-equivariant where the $W$ action
on $(G/B)^T=W$ is given by $w\cdot (x)_{x\in W}:=(xw)$. Moreover, the
$W$-action on $(G/B)^T$ induces a map on $K_{T_{comp}}^0((G/B)^{T})$
given by $w\cdot (a_x)_{x\in W}:=(a_{xw})_{x\in W}$. Recall that this
action of $W$ on $K_{T_{comp}}^0(G/B)^{T}$ restricts to the above
action on $K_{T_{comp}}^0(G/B)$ (see \cite[Section 1.4]{ml},
\cite[Definition 3.11]{kk}) which we call as the {\it star}
action. Moreover, the map $\Psi$ defined earlier (see \eqref{characteristic}) is $W$-equivariant
for the star action of $W$ on $K_{T_{comp}}^0(G/B)$ and the $W$ action
on
$R(T_{comp})\otimes_{R(G_{comp})}
K^0_{G_{comp}}(G/B)=R(T_{comp})\otimes_{R(G_{comp})} R(T_{comp})$
given by $w\cdot \sum c_i\otimes d_i:=\sum c_i\otimes w\cdot d_i$ (see
\cite{ml}). In particular, the star action of each $w\in W$ on
$K_{T_{comp}}^0(G/B)$ is a $R(T_{comp})$-algebra automorphism (see
\cite{kk}), unlike the dot action of $w\in W$ which is a twisted
$R(T_{comp})$-algebra automorphism in the sense that
$w\cdot (cf)=w(c)w(f)$ for $c\in R(T_{comp})$ and
$f\in K_{T_{comp}}^0(G/B)$. See \cite[Theorem 3.1, Theorem 3.2]{tym}
for the analogous statements for the equivariant cohomology of $G/B$.

\subsubsection{Some additional properties when $\pi_1(G_{comp})$ is free}\label{free}

In this section we shall assume in addition that $\pi_1(G_{comp})$ is free.

The following lemma is the analogue for $G_{comp}$-equivariant topological $K$-ring of a $T$-cellular $G$-variety of \cite[Lemma 1.6]{u1}.

\blem
The ring $K_{G_{comp}}^0(X)$ is a free $R(G_{comp})$-module of rank $m=|X^{T}|$.
\elem

\begin{proof}
  Recall that $R(G_{comp})\cong R(T_{comp})^{W}$ and $R(T_{comp})$ is
  a free $R(G_{comp})$-module of rank $|W|$. Furthermore, since
  $K^0_{T_{comp}}(X)$ is a free $R(T_{comp})$-module of rank
  $|X^{T}|=m$ it follows that $K^0_{T_{comp}}(X)$ is a free
  $R(G_{comp})$-module. Furthermore, from the isomorphism $\Psi$ (see \eqref{characteristic}) we
  get that $K_{G_{comp}}^0(X)$ is a direct summand of
  $K^0_{T_{comp}}(X)$ as an $R(G_{comp})$-module. Thus it is a
  projective and is hence a free $R(G_{comp})$-module, since
  $R(G_{comp})$ is a tensor product of a polynomial ring and a Laurent
  polynomial ring and hence is a regular ring (see \cite{pit},
  \cite{gub}).

  By the isomorphism $\Psi$ it follows that the rank of
  $K_{G_{comp}}^0(X)$ as a free $R(G_{comp})$-module is $m$.
  \end{proof}

  Let $W^I$ denote the set of minimal length coset representatives of
  the parabolic subgroup $W_{I}$ for every $I\subseteq \Delta$.  Note
  that $J\subseteq I$ implies that
  $W^{\Delta\setminus J}\subseteq W^{\Delta\setminus I}$.  Let
  $\displaystyle C^{I}:=W^{\Delta\setminus I}\setminus \Big (
  \bigcup_{J\subsetneq I} W^{\Delta\setminus J} \Big)$

  Recall that in \cite{st} Steinberg has defined a basis
  $\{f_v^{I}:v\in W^{I}\}$ for $R(G_{comp})=R(T_{comp})^{W}$ as
  $R(T_{comp})$-module and $R(G_{comp})^{W_{I}}$ as
  $R(T_{comp})$-module for $I\subseteq \Delta$ (also see \cite[Section
  1]{u2} for more details).

  Whenever $v\in C^{I}$ we denote $f_{v}^{W_{\Delta\setminus I}}$
  simply by $f_v$. By \cite[Lemma 1.10]{u1},
  $\{f_v:v\in \sqcup_{J\subseteq I} C^{J}=W^{\Delta\setminus I}\}$ is
  a $R(G_{comp})=R(T_{comp})^{W}$-basis of
  $R(T_{comp})^{W_{\Delta\setminus I}}$ for every $I\subseteq \Delta$. We let

  \be\label{dec1} R(T_{comp})_{I}:=\bigoplus_{v\in C^{I}}
  R(G_{comp})\cdot f_v.\ee

  In $R(T_{comp})$ we write

  \be\label{mult1} f_v\cdot f_{v'}=\sum _{J\subseteq I\cup I'}\sum _{w\in
    C^{J}}a^{w}_{v,v'} \cdot f_{w}\ee for $a^w_{v,v'}\in R(G_{comp})$,
  where $v\in C^{I}$ and $v'\in C^{I'}$ and $w\in C^{J}$,
  $J\subseteq I\cup I'$.

\subsection{Relation with ordinary topological $K$-theory}\label{relord}

Recall from \cite[Corollary 4.4]{u4} that $T$-cellular varieties are {\em
  weakly equivariantly formal in $K$-theory} with respect to
$T_{comp}$-action, in the sense of \cite{hl}. More precisely, we have the following theorem:

\bcor\label{tequivformal}(\cite[Corollary 4.4]{u4}) We have the isomorphism \be\label{teqf}\mathbb{Z}\otimes_{R(T_{comp})}{K}^0_{T_{comp}}(X) \cong {K}^0(X).\ee
Here the map from $K_{T_{comp}}(X)\lra K(X)$ is the canonical
forgetful map, $\mathbb{Z}\lra K^0(X)$ is the map which takes an
integer $n$ to the class of the trivial vector bundle of rank $n$ and
$\mathbb{Z}$ is $R(T_{comp})$-module under the augmentation map. \ecor

Moreover, from \cite[Lemma 4.4]{hl} we have the following:

\bth\label{haradalandweber} (\cite[Lemma 4.4]{hl}) Let $G_{comp}$ be a compact connected Lie group with $\pi_1(G_{comp})$ torsion-free, and let $M$ be a
compact $G$-space. Then $M$ is weakly equivariantly formal with respect to $G_{comp}$ if and only if it is weakly equivariantly formal with respect to a maximal torus $T_{comp} \subset G_{comp}$.
\eeth

We have the following theorem which follows immediately from Theorem \ref{tequivformal} and Theorem \ref{haradalandweber}.

\bth\label{GEQF} A $T$-cellular $G$-variety with $\pi_1(G_{comp})$ free is weakly equivariantly formal in $K$-theory with respect to the action of $G_{comp}$ in the sense of \cite{hl}. More precisely, we have the following isomorphism:
\be\label{geqf}\mathbb{Z}\otimes_{R(G_{comp})}{K}^0_{G_{comp}}(X) \cong {K}^0(X).\ee
\eeth

\section{Cellular Toric varieties}\label{cellulartv}

In this section we shall recall from \cite[Theorem 3.1]{u3} a
combinatorial characterization on the fan $\Sigma$ for the toric
variety $X(\Sigma)$ to be $T$-cellular (see Definition
\ref{cellular}). We call a fan $\Sigma$ satisfying the necessary and
sufficient combinatorial conditions a {\em cellular fan} (see
\cite[Definition 3.2]{u3}).

We begin by fixing some notations and conventions.

Let $X=X(\Sigma)$ be the toric variety associated to a  fan
$\Sigma$ in the lattice $N\simeq \bz^n$. Let $M:=Hom(N,\bz)$ be the
dual lattice of characters of $T$. Let $\{v_1,\ldots,v_d\}$ denote the
set of primitive vectors along the edges
$\Sigma(1):=\{\rho_1,\ldots,\rho_d\}$. Let $V(\gamma)$ denote the
orbit closure in $X$ of the $T$-orbit $O_{\gamma}$ corresponding to
the cone $\gamma\in\Sigma$. Let $S_{\sigma}=\sigma^{\vee}\cap M$ and
$U_{\sigma}:=Hom_{sg}(S_{\sigma},\mathbb{C})$ denote the $T$-stable
open affine subvariety corresponding to a cone $\sigma\in \Sigma$.

We further assume that all the maximal cones in $\Sigma$ are
$n$-dimensional, in other words $\Sigma$ is {\it pure}. 
The $T$-fixed locus in $X$ consists of the set of $T$-fixed points
\be\label{fp}\{x_{1}, x_{2}\ldots,x_{m}\}\ee
corresponding to the set of maximal dimensional cones
\be\label{maxc}\Sigma(n):=\{\sigma_1,\sigma_2,\ldots,\sigma_m\}.\ee

Choose a {\it generic} one-parameter subgroup $\lambda_v\in X_{*}(T)$
corresponding to a $v\in N$ which is outside the hyperplanes spanned
by the $(n-1)$-dimensional cones, so that (\ref{fp}) is the set of
fixed points of $\lambda_v$ (see \cite[\S3.1]{Br2}).  For each
$x_{i}$, we have the plus stratum
$X_{\lambda_{v},x_i}=Y_i$ (see \eqref{ps}).

Consider a face $\gamma$ of $\sigma_i$ satisfying the property that
the image of $v$ in $N_{\br}/\br\gamma$ is in the relative interior of
$\sigma_i/\br\gamma$. Since the set of such faces is closed under
intersections, we can choose a minimal such face of $\sigma_i$ which
we denote by $\tau_i$. We have \be\label{pu}
Y_i=\bigcup_{\tau_i\subseteq \gamma\subseteq \sigma_i}O_{\gamma}.\ee
Moreover, if we choose a generic vector $v\in|\Sigma|$, then
\[\displaystyle{X=\bigcup_{i=1}^m Y_i}\] where $X_{\lambda_{v},x_i}=Y_i$,
$1\leq i\leq m$ are the cells of the {\it Bialynicki-Birula cellular
  decomposition} of the toric variety $X$ corresponding to
the one-parameter subgroup $\lambda_v$ (see \cite[Lemma 2.10]{hs}).

Furthermore, the Bialynicki-Birula decomposition for $X$ corresponding
to $\lambda_v$ is {\it filtrable} if
\be\label{clstar}\bar{Y_{i}}\subseteq \bigcup_{j\geq i} Y_{j}\ee for
every $1\leq i\leq m$.  It can be seen that \eqref{clstar} is
equivalent to the following combinatorial condition in $\Sigma$:
\be\label{star} \tau_i\subseteq \sigma_j~~\mbox{implies}~~i\leq j.\ee
Now, from \eqref{pu} it follows that
$Y_i=V(\tau_i) \cap U_{\sigma_i}$. Let $k_i:=\dim(\tau_i)$ for
$1\leq i\leq m$. Thus $Y_i$ is a $T$-stable affine open set in the
toric variety $V(\tau_i)$, corresponding to the maximal dimensional
cone $\bar{\sigma_i}$ in the fan
$(\mbox{star}(\tau_i), N(\tau_i)=N/N_{\tau_i})$ (see \cite[Chapter
3]{f}). Thus $Y_i$ is isomorphic to the complex affine space
$\mathbb{C}^{n-k_i}$ if and only if the $(n-k_i)$-dimensional cone
$\bar{\sigma_i}:=\sigma_i/N_{\tau_i}$ is a smooth cone in
$(\mbox{star}(\tau_i), N(\tau_i))$ for $1\leq i\leq m$.

It follows by \eqref{clstar} that under the above conditions on
$\Sigma$, namely \eqref{star} and that $\bar{\sigma_i}$ is a smooth cone in
$\mbox{star}(\tau_i)$, $\displaystyle Z_i:=\bigcup_{j\geq i} Y_j$ are
$T$-stable closed subvarieties, which form a chain
\be\label{strata}X=Z_1\supseteq Z_2\supseteq \cdots \supseteq
Z_m=V(\tau_m)\ee such that
$Z_{i}\setminus Z_{i+1}=Y_i\simeq \mathbb{C}^{n-k_i}$ for
$1\leq i\leq m$.

We summarize the above discussion in the following theorem which gives
a combinatorial characterization on $\Sigma$ for $X(\Sigma)$ to be a
$T$-cellular toric variety (see Definition \ref{cellular}).

Fix a generic $v\in |\Sigma|\cap N$ such that $\lambda_v$ is a generic
one-parameter subgroup of $T$.

\bth\label{combcell} The toric variety $X(\Sigma)$ is $T$-{\it
  cellular} if and only if the following combinatorial conditions hold
in $\Sigma$:

\noindent
(i) $\Sigma$ admits an ordering of the
maximal dimensional cones
\be\label{ordering}\sigma_1<\sigma_2<\cdots<\sigma_m\ee such that the
distinguished faces $\tau_i\subseteq \sigma_i$ defined above satisfy
the following property: \be\label{star1} \tau_i\subseteq
\sigma_j~~\mbox{implies}~~i\leq j 
\ee
\noindent
(ii) $\bar{\sigma_i}:=\sigma_i/N_{\tau_i}$ is a smooth cone in the fan
$(\mbox{star}(\tau_i), N(\tau_i))$ for $1\leq i\leq m$.
\eeth

The following technical lemma shows that for a toric variety $X(\Sigma)$ we can
vary the generic one-parameter subgroup suitably to get the same
Bialynicki-Birula cell decomposition. This shall be used in Section
\ref{celltoremb}.

\blem \label{lem:celluparam} Let $X=X(\Sigma)$ be a toric $T$-variety
associated to a fan $\Sigma$ in $N=X_*(T)\simeq \bz^n$. Let
$\sigma\in \Sigma(n)$ be a maximal dimensional cone. Let
$\mu_{1,\sigma},\ldots,\mu_{q_\sigma,\sigma} \in X^*(T)$ denote a
finite number of generators of the semigroup $\sigma^\vee$ :
  \[\sigma^\vee = \bz_{\ge 0}\mu_{1,\sigma}+\cdots+\bz_{\ge 0}\mu_{q_\sigma,\sigma}\, .\]

  Let $\lambda\in N$ be a one-parameter subgroup which is generic for
  $X$ (i.e. $X^\lambda=X^T$). If $\lambda'\in N$ is such that
  \be\label{eq:condgen}~\forall~ i~,\, ~|\langle \mu_{i,\sigma} ,
  \lambda-\lambda'\rangle|<|\langle\mu_{i,\sigma},\lambda\rangle|\ee
  then $\lambda'$ is also a generic for $X$ and the Bialynicki-Birula
  cells of $X$ corresponding to $\lambda$ and $\lambda'$ are the same:
  \[X_{\lambda,x_{\sigma}}= X_{\lambda',x_{\sigma}},\] where
  $x_\sigma$ is the $T$-fixed point associated to the cone $\sigma$.
  In particular, if $X$ is $T$-cellular with respect to $\lambda$ then
  it is $T$-cellular with respect to $\lambda'$ and vice versa. \elem

\begin{proof} If $\gamma \subset\sigma$ is a face, let \[\mathrm{Int}(\sigma/\br\gamma)\]
  be the relative interior of $\sigma/\br\gamma$ in the quotient space $N_\br/\br \gamma$.

  Then one can describe the relative interior as :
  \[\mathrm{Int}(\sigma/\br\gamma) = \{\bar{x}\, : \,  x \in N_\br,\, \forall ~\mu \in \sigma^\vee\setminus \gamma^\perp,\, \langle \mu, x\rangle >0\}\]
  \[=\{\bar{x}\, : \,  x \in N_\br,\, \forall~ 1\le i\le q_\sigma,\, \mu_{i,\sigma}\not\in \gamma^\perp \Rightarrow  \langle \mu_{i,\sigma}, x\rangle >0\}\, .\]

  Because of the condition (\ref{eq:condgen}),   the real numbers
  \[ \langle \mu_{i,\sigma},\lambda\rangle,\, \langle \mu_{i,\sigma}, \lambda'\rangle \]
  have the same sign for all $1\le i \le q_\sigma$. Therefore \[\lambda\in \mathrm{Int}(\sigma/\br\gamma)\Leftrightarrow \lambda' \in \mathrm{Int}(\sigma/\br\gamma)\]
  for each face $\gamma\subset\sigma$.

  In particular, the minimal face $\gamma$ is such that
  $\lambda \in \mathrm{Int}(\sigma/\br\gamma)$ is also the minimal
  face such that $\lambda'\in \mathrm{\sigma/\br\gamma}$. Therefore one has
  \[X_{\lambda,x_\sigma}=X_{\lambda',x_\sigma}\, .\]
\end{proof}

\brem\label{tau_1} We shall assume without loss of generality that $v$
belongs to the relative interior of $\sigma_1$ so that
$\tau_1=\{0\}$. In particular, this implies that $\sigma_1$ is a
smooth cone in $\Sigma$ so that
$Y_1=U_{\sigma_1}\cap V(\tau_1)\simeq \mathbb{C}^n$ is smooth in
$X=V(\tau_1)$.  \erem

\subsection{\boldmath $T_{comp}$-equivariant topological $K$-ring of $X(\Sigma)$}  
From the description of the $T$-fixed points and the $T$-stable curves
in $X=X(\Sigma)$ (see \cite{f}) it follows that the set $\mathcal{A}$
defined in Section \ref{GKMsection}, is the set of all
$(a_i)_{1\leq i\leq m}\in R(T_{comp})^m$ such that
$a_i\equiv a_j \pmod{(1-e^{\chi})}$ whenever the maximal dimensional
cones $\sigma_i$ and $\sigma_j$ share a wall $\sigma_i\cap \sigma_j$
in $\Sigma$ and $\chi\in (\sigma_i\cap \sigma_j)^{\perp}\cap M$. We
further recall from \cite[proof of Theorem 5.2]{u3}) that a complete
$T$-cellular toric variety $X=X(\Sigma)$ satisfies Assumption
\ref{as}.

Thus by Theorem \ref{main1}, we get the following GKM type description
for $K^0_{T_{comp}}(X)$ as a
$K^0_{T_{comp}}(pt)=R(T_{comp})$-subalgebra of $R(T_{comp})^m$. For a
direct proof of this result we refer to \cite{u3}.

\bth\label{tveqktdes1}(\cite[Theorem 5.2]{u3}) Let $X=X(\Sigma)$ be a
complete $T$-cellular toric variety. The ring $K^0_{T_{comp}}(X)$ is
isomorphic to $\mathcal{A}$ as an $R(T_{comp})$-subalgebra of
$R(T_{comp})^m$.  \eeth

\subsubsection{The ring of piecewise Laurent polynomial functions on
  $\Sigma$}\label{piecewiselaurent}

Let $\Sigma$ be a fan in $N\simeq \mathbb{Z}^n$.  For
$\sigma\in \Sigma$ let $T_{\sigma}\subseteq T$ denote the stabilizer
along $O_{\sigma}$. Since $X^*(T_{\sigma})=M/M\cap \sigma^{\perp}$ we
can identify $R(T_{\sigma})$ with the ring
$\mathbb{Z}[M/M\cap \sigma^{\perp}]$ of Laurent polynomial functions
on $\sigma$.  Furthermore, whenever $\sigma$ is a face of
$\sigma'\in\Sigma$, noted $\sigma\preceq \sigma'$, we have a natural homomorphism
$\varphi_{\sigma,\sigma'}:\mathbb{Z}[M/\sigma'^{\perp}\cap
M]=R(T_{\sigma'})\ra \mathbb{Z}[M/\sigma^{\perp}\cap
M]=R(T_{\sigma})$, given by restriction of Laurent polynomial
functions on $\sigma'$ to $\sigma$. Let
\be\label{plp}PLP(\Sigma):=\{(f_{\sigma}) \in \prod_{\sigma\in \Sigma}
{\mathbb{Z}[M/M\cap \sigma^{\perp}]}\mid\varphi_{\sigma,\sigma'}
(f_{\sigma'})= f_{\sigma} ~\mbox{whenever}~\sigma\preceq \sigma' \in
\Sigma\}.\ee Then $PLP(\Sigma)$ is a ring under pointwise addition and
multiplication and is called the ring of piecewise Laurent polynomial
functions on $\Sigma$. Moreover, we have a canonical map
$R(T_{comp})=\mathbb{Z}[M]\ra PLP(\Sigma)$ which sends $f$ to the
constant tuple $(f)_{\sigma\in \Sigma}$. This gives $PLP(\Sigma)$ the
structure of $R(T_{comp})$-algebra. The following result was proved in
\cite{u3}.

\bth\label{tveqkthdes2} (see \cite[Theorem 5.6]{u3}) Let $X=X(\Sigma)$
be a complete $T$-cellular toric variety.  The ring
$K^0_{T_{comp}}(X(\Sigma))$ is isomorphic to $PLP(\Sigma)$ as an
$R(T_{comp})$-algebra.  \eeth

\section{Cellular toroidal embeddings}\label{celltoremb}

In this section we shall describe cellular toroïdal embeddings and show
that they satisfy Assumption \ref{as}.

We isolate below the definition of toroidal $G$-embeddings given in
Section \ref{Introduction}.

\bdefe\label{toroidal}(\cite[Definition 6.2.2]{bk})  A $G$-embedding
$\mathbb{X}$ is {\em toroidal} if the quotient map $p:G\lra G_{ad}$ extends
to a unique $G\times G$-equivariant morphism 
$p:\mathbb{X}\lra \bar{G_{ad}}$ where $\bar{G_{ad}}$ denotes the wonderful
compactification of $G_{ad}$ defined by De Concini and Procesi in
\cite{dp}. Here, the $G\times G$-action on $\bar{G_{ad}}$ is via its
projection to $G_{ad}\times G_{ad}$.\edefe

We shall follow the notation in Section \ref{Introduction} together with the following.

Let $T_{ad}:=p(T)$ be the quotient of $T$ by its center $C$ and let
$\bar{T_{ad}}$ be the closure of $T_{ad}$ in $\bar{G_{ad}}$. Recall
from \cite[Lemma 6.1.6]{bk} that $\bar{T_{ad}}$ is a nonsingular
$T_{ad}$-toric variety whose fan $\mathcal{F}_{ad}$ in
$X_*(T_{ad})\otimes \mathbb{R}$ consists of the Weyl chambers and
their faces.

The action of $W$ on $X_*(T_{ad})\otimes \mathbb{R}$ preserves the fan $\mathcal{F}_{ad}$. 
Recall that the $\mbox{diag}(W)$ action on $T_{ad}$ extends to
$\bar{T_{ad}}$. This corresponds to the action of $W$ on the fan
$\mathcal{F}_{ad}$. Let
$\mathcal{F}_{ad}^+$ denote the fan associated to the positive Weyl
chamber. The toric variety $\bar{T_{ad}}^+$ associated to
$\mathcal{F}_{ad}^+$ is isomorphic to the affine space $\mathbb{A}^r$
where $T_{ad}$ acts on $\mathbb{A}^r$ by the embedding
$t\mapsto (\alpha_1(t),\ldots, \alpha_r(t))$.  Furthermore,
$\mathcal{F}_{ad}=W\mathcal{F}_{ad}^+$ so that
$\bar{T_{ad}}=\mbox{diag}(W)\bar{T_{ad}}^+$.

Let $\mathbb{X}$ be a toroidal embedding of $G$ (see Definition
\ref{toroidal}).  Let $X:=\bar{(T\times T)\cdot x}$ be the closure of
$T$ in $\mathbb{X}$. Then $X$ is a $T$-toric variety satisfying
$X=p^{-1}(\bar{T_{ad}})$. Moreover, if $X^+:=p^{-1}(\bar{T_{ad}}^+)$
then $X=\mbox{diag}(W) X^+$. Then $X^+$ is the $T$-toric variety whose
fan $\cf^+$ corresponds to a subdivision of the positive Weyl chamber
in $X_*(T)\otimes \mathbb{R}$. In addition, the fan $\cf$ of $X$
satisfies \be\label{permutation} \cf = W \cf ^+.\ee Furthermore, $p$
restricts to a morphism from $X$ to $\bar{T_{ad}}$ which is
equivariant with respect to the actions of $T$ and $W$. Moreover, the
restriction of $p$ to $X^+$ defines a $T$-equivariant morphism from
$X^+\ra \bar{T_{ad}}^+$ and hence corresponds to a morphism of fans
\be (\mathcal{F}^+, X_*(T))\ra (\mathcal{F}_{ad}^+,X_*(T_{ad})).\ee We
refer to \cite[Proposition 6.2.3, Proposition 6.2.4]{bk} for a precise
characterization of the geometry of $\mathbb{X}$ in relation to the
geometry of $X^+$ and hence that of $X$. In particular, by
\cite[Proposition 6.2.4 (iii)]{bk} the toroidal embeddings of $G$ are
classified by the fans in $X_*(T)\otimes \mathbb{R}$ with support in
the positive Weyl chamber.

The closed $G\times G$ orbits of $\mathbb{X}$ are each isomorphic to $G/B^{-}\times G/B$ with distinguished points $x_{\sigma}$ corresponding to
$\sigma\in \cf^+(l)$. Moreover, the $T\times T$-fixed points of
$\mathbb{X}$ all lie in the closed orbits and consist of
$(w_1,w_2)\cdot x_{\sigma}$ where $(w_1,w_2)\in W\times W$ and
$x_{\sigma}$ for $\sigma\in \cf^{+}(l)$ are the $T$-fixed points of
$X^+$. In particular, $|\mathbb{X}^{T\times T}|=m\cdot |W|^2$ where
$m=|\cf^{+}(l)|$.

As in Section \ref{Introduction}, let $\mathbb{X}$ be a toroidal
embedding of $G$ and let $X=\bar{T\times T\cdot x}$ the toric
$T$-variety which characterizes $\mathbb{X}$.

For any $\sigma \in \mathcal{F}(l)$ let
$\mu_{1,\sigma},\ldots,\mu_{q_\sigma,\sigma} \in X^*(T)$ be a finite
number of generators of the semigroup $\sigma^\vee$ :
  \[\sigma^\vee = \bz_{\ge 0}\mu_{1,\sigma}+\cdots +\bz_{\ge 0}\mu_{q_\sigma,\sigma}\, .\]

  In the following proposition we show that if all the
  Bialynicki-Birula cells of the $T\times T$-variety $\mathbb{X}$ are
  smooth for a suitable choice of a one-parameter subgroup of
  $T\times T$ then for the toric $T$-variety $X=X(\mathcal{F})$, we
  can find a suitable one-parameter subgroup of $T$ for which the
  Bialynicki-Birula cells are all smooth and vice-versa.

More precisely one has the following statements.

\bpropo\label{equivalentcellularcondition} \begin{enumerate}
\item If for the toroidal embedding $\mathbb{X}$ the Bialynicki Birula
  cells corresponding to a one parameter subgroup
  $(\nu_1,\nu_2)\in X_*(T)\times X_*(T)$ are all smooth then for the
  toric $T$-variety $X$ the Bialynicki-Birula cells corresponding to
  the one parameter subgroup
  $w_1\cdot \nu_1-w_2\cdot \nu_2 \in X_*(T)$ are all smooth, where
  $w_1,w_2\in W$ are such that $w_1\cdot \nu_1$ is dominant and
  $w_2\cdot \nu_2$ is anti-dominant.
\item If for the toric $T$-variety $X$ the Bialynicki-Birula cells
  corresponding to a one parameter subgroup $\nu_0\in X_*(T)$ are all
  smooth, then for the $T\times T$-variety $\mathbb{X}$ the
  Bialynicki-Birula cells corresponding to the one parameter subgroup
  $(\nu_1,\nu_2)\in X_*(T)\times X_*(T)$ are all smooth, where
  $\nu_2\in X_*(T)$ is any regular one parameter subgroup and
  $\nu_1=N\nu_0$ with $N$ an integer such
  that\begin{equation}\label{eq:cond_generic} N>
    max\left\{\frac{|\langle \mu_{i,\sigma},w\cdot \nu_2\rangle |}{|\langle
        \mu_{i,\sigma},\nu_0\rangle |} \, : \, \sigma \in
      \mathcal{F}(l),\, 1\le i \le q_\sigma \right\}
    \end{equation}
\end{enumerate}\epropo

\begin{proof}

  {\it Proof of (1):} Suppose that there are $\nu_1,\nu_2\in X_*(T)$
  one-parameter subgroups such that
  $\nu:=(\nu_1,\nu_2)\in X_*(T\times T)$ is generic i.e.,
\[\mathbb{X}^{\nu}=\mathbb{X}^{T \times T}=\{(w,w')\cdot x_\sigma\mid \sigma\in \cf^+(l)~ (w,w')\in W\times W\}\] the set of $T\times T$-fixed
points of $\mathbb{X}$. Furthermore, for each $\sigma\in \cf_+(l)$ and
$(w,w')\in W\times W$, the Bialynicki-Birula cell
\[\mathbb{X}_{\nu,(w,w')\cdot x_{\sigma}}:=\{ x\in \mathbb{X}\mid \lim_{t\ra 0}\nu(t)x ~
  \mbox{exists and is equal to}~(w,w')\cdot x_{\sigma}\},\] associated
to the one-parameter subgroup $\nu$ and the fixed point
$(w,w')\cdot x_{\sigma}$ is smooth. Let $w_1,w_2\in W$.  A quick calculation
show that
\[(w_1,w_2).\mathbb{X}_{\nu,x_{\sigma}}=\mathbb{X}_{(w_1\cdot \nu_1,w_2\cdot
    \nu_2), (w_1,w_2)\cdot x_{\sigma}}.\]

In particular, it follows that the Bialynicki-Birula cells of
$\mathbb{X}$ corresponding to the one-parameter subgroup
$(w_1\cdot \nu_1,w_2\cdot \nu_2)$ are also smooth.

So by choosing convenient $w_1, w_2\in W$, we can assume that $\nu_1$
is dominant and $\nu_2$ is anti-dominant\footnote{\ie
  $\forall \alpha\in \Phi^+,\, \langle \alpha,\nu_1\rangle >0,\,
  \langle\alpha,\nu_2\rangle <0$.}

Let $(\bar{G_{ad}})_0$ denote the $B_{ad}\times B^-_{ad}$-big cell of
the wonderful compactification $\bar{G_{ad}}$ that is the
Bialynicki-Birula open cell associated to the one-parameter subgroup
$\nu$ of $T\times T$.

Put $\mathbb{X}_0:=p^{-1}( (\bar{G_{ad}})_0),\, X_0:=X\cap \mathbb{X}_0$.

Then $X_0=X^+$ is the toric variety associated to the fan
$\mathcal{F}^+$, which is the intersection of $\mathcal{F}$ with the
positive Weyl chamber of $X_*(T)_{\mathbb{R}}$
  
According to \cite{bk}[Prop. 6.2.3.(i)], the map
  \begin{eqnarray}\label{eq:grocell}
    \varphi\, : \, U \times U^-\times X_0\to \mathbb{X}_0,\, (u,v,z)\mapsto (u,v).z\end{eqnarray}
  is an isomorphism of $T\times T$-varieties, where $(t_1,t_2)\in T\times T$ acts over $U$ by $t_1$, over $U^-$ by $t_2$ and by $t_1t_2^{-1}$ over $X_0$. But $X_0=X(\mathcal{F}^+)$ is open in $X=X(\mathcal{F})$. If  $\sigma$ is a maximal cone of $\mathcal{F}^+$ then the corresponding point $x_\sigma$ is in $X_0$. Thus the Bialynicki-Birula cell ${X}_{\nu_1-\nu_2, x_\sigma}$ of $X$,  is contained in $X_0$. Likewise, since $\mathbb{X}_0$ is an open subset of $\mathbb{X}$ that contains $x_\sigma$, $\mathbb{X}_{\nu,x_\sigma} \subseteq \mathbb{X}_0$. Moreover, $\mathbb{X}_{\nu,x_\sigma}$ is $B \times B^-$-stable (because $\nu_1$ is dominant and $\nu_2$ is anti-dominant). So, one has 
  \[\varphi(U\times U^-\times {X}_{\nu_1-\nu_2,x_\sigma})=\mathbb{X}_{\nu,x_\sigma}\]
  \[\implies U\times U^-\times {X}_{\nu_1-\nu_2,x_\sigma}\cong \mathbb{X}_{\nu,x_\sigma}\]

  In particular, if $\mathbb{X}_{\nu,x_\sigma}$ is smooth then so is $X_{\nu_1-\nu_2,x_\sigma}$. 

  So if the Bialynicki-Birula cells of $\mathbb{X}$ for the
  one-parameter subgroup $(\nu_1,\nu_2)$ of $T\times T$ are smooth
  implies that the Bialynicki-Birula cells of $X$ for the
  one-parameter subgroup $\nu_1-\nu_2$ of $T$ are all smooth.

  \nopagebreak[4] {\it Proof of (2):} Suppose that for the toric
  $T$-variety $X$ all the Bialynicki-Birula cells corresponding to a
  one-parameter subgroup $\nu_0\in X_*(T)$ are smooth. For every
  $w \in W$, and each $T$-fixed point $x\in X$ we have,
\begin{eqnarray}\label{eq:Xwnu}
{X}_{\nu,x} \simeq {X}_{w\cdot \nu, w\cdot x}\, .
\end{eqnarray}

It follows that for every $w \in W$, the Bialynicki-Birula cells of
$X$ corresponding to $w\cdot \nu$ are also smooth. Since one can take
$w\cdot \nu$ instead of $\nu$, one can suppose that $\nu$ is
dominant. It is understood that $\nu$ is a regular dominant
one-parameter subgroup since we assume that ${X}^\nu={X}^{T}$.

Let $\sigma\in \mathcal{F}^+(l)$ be a maximal cone of dimension $l$ in
the fan $\mathcal{F}^+$. Let $n_1,n_2\in N(T)$ and let
$w_1=n_1T,w_2=n_2T$ be the corresponding elements in $W$. The
isomorphism $\varphi$ of (\ref{eq:grocell}) induces an
isomorphism \begin{equation}\label{eq:cellaff} n_1U \times
  n_2U^-\times X_0\to (n_1,n_2)\cdot \mathbb{X}_0,\, (u,v,z)\mapsto
  (u,v)\cdot z\, .\end{equation}

Let $\nu_1,\nu_2\in X_*(T)$ such that $\nu_1$ is regular dominant and
$\nu_2$ regular anti-dominant.
  
Since the  $T\times T$-fixed point
\[(n_1,n_2)\cdot x_\sigma \in (n_1,n_2)\cdot \mathbb{X}_0,\] since
$(n_1,n_2)\cdot \mathbb{X}_0$ is open in $\mathbb{X}$, one has the
inclusion
\[\mathbb{X}_{(\nu_1,\nu_2),(n_1,n_2)\cdot x_\sigma
  }\subset(n_1,n_2)\cdot \mathbb{X}_0\, .\]

  But, if $u\in U,\, v\in U^-$, if $y \in X_0$, then
  \[(n_1u, n_2v)\cdot y\in X_{(\nu_1, \nu_2),(n_1,n_2)\cdot x_\sigma
    }\]\[\Leftrightarrow u\in U\cap w_1^{-1} Uw_1,\, v \in U^-\cap
    w_2^{-1} U^-w_2,\, y\in {X}_{w_1^{-1} \cdot \nu_1-w_2^{-1} \cdot
      \nu_2, x_\sigma}\, .\]

  Then one can deduce the isomorphism
  \[\mathbb{X}_{(\nu_1,\nu_2),(w_1,w_2)\cdot x_\sigma }\simeq U\cap
    w_1^{-1} U w_1 \times U^-\cap w_2^{-1}U^-w_2 \times {X}_{w_1^{-1}
      \cdot \nu_1-w_2^{-1} \cdot \nu_2, x_\sigma}\, .\]

  It is therefore sufficient to show that, for a good choice of
  $\nu_1,\nu_2$, the cells of the toric variety $X$
  \[{X}_{w_1^{-1} \cdot \nu_1-w_2^{-1} \cdot\nu_2,x_\sigma}\] are
  smooth for all $w_1,w_2\in W$ and every maximal dimensional cone
  $\sigma$ in $\mathcal{F}^+(l)$.

  Since for all $w\in W$, the isomorphism $X\to (w,w)X=X$ which maps
  $x\mapsto (w,w)\cdot x$ induces an isomorphism of the
  Bialynicki-Birula cells :
  \[X_{\nu,x_{\sigma}}\simeq X_{w\cdot \nu,w\cdot x_{\sigma}}\] for
  all one-parameter subgroup $\nu \in X_*(T)$ and all
  $\sigma\in \mathcal{F}(l)$, it is sufficient to prove that, for a
  convenient choice of $\nu_1,\nu_2$ the cells
  \begin{eqnarray}\label{eq:lisses}
    X_{\nu_1-w_1w_2^{-1}\nu_2, x_\sigma}
\end{eqnarray}
are smooth for all $w_1,w_2\in W$ and all $\sigma\in \mathcal{F}(l)$.

We suppose that the dominant one-parameter subgroup $\nu_0$ is generic
for the toric variety $X$.

  It means that $X^{\nu_0}=X^T$ that is :
  \[\forall ~\sigma\in \mathcal{F}(l),\, ~\forall ~ 1\le i \le q_\sigma,\,
    \langle \mu_{i,\sigma},\nu_0\rangle \neq 0\, .\]

  {\it Remark.} As $\mathcal{F}$ is a subdivision of the Weyl
  chambers, the one-parameter subgroup $\nu_0$ is regular.
  
  Let $\nu_2\in X_*(T)$ be any regular anti-dominant one-parameter
  subgroup and $N$ be an integer as in (\ref{eq:cond_generic}).
  
  Let $\nu_1=N\nu_0$. Then $\nu=(\nu_1,\nu_2)$ satisfies the condition
  (\ref{eq:lisses}).

  The choice of $N$ implies that
  \[\forall~ w \in W,\,~\forall ~\sigma\in \mathcal{F}(l),\, ~ \forall~ 1\le
    i\le
    q_\sigma~,\,\]\[ |\langle \mu_{i,\sigma},w\cdot \nu_2\rangle|<
    |\langle \mu_{i,\sigma},N\nu_0\rangle | .\]

  In particular, according to the Lemma \ref{lem:celluparam}, if
  $\nu_0$ is generic for the toric $T$-variety $X$,
  $\nu_0-w\cdot \nu_2$ is also generic for every $w \in W$, and the
  Bialynicki-Birula cells
  $X_{\nu_0,x_\sigma}=X_{\nu_1-w\cdot \nu_2,x_\sigma}$ are the same
  for all $w\in W$.

  As the Bialynicki-Birula cells of the toric $T$-variety $X$
  corresponding to $\nu_0$ are smooth the same is true for the
  Bialynicki-Birula cells associated to the one-parameter subgroup
  $w_1\nu_1-w_2\nu_2$ for all $w_1,w_2\in W$.

  As $\nu_1$ is regular dominant and $\nu_2$ regular antidominant, we
  have shown that the Bialynicki-Birula cells of the
  $T\times T$-variety $\mathbb{X}$ corresponding to the generic
  one-parameter subgroup $\nu=(\nu_1,\nu_2)=(N\nu_0,\nu_2)$ of
  $T\times T$ are smooth.

\end{proof}

 Suppose $\cf^+(l)=\{\sigma_1,\ldots, \sigma_m\}$. Let $x_i$ denotes
  the $T$-fixed point of $X^+$ associated to $\sigma_j\in\cf^+(l)$ for
  $1\leq j\leq m$. Then we have \[(X^+)^{T}=\{x_1,\ldots, x_m\},\]
  \[\mathbb{X}^{T}=\{(w_1,w_2)\cdot x_i\mid w_1,w_2\in W, 1\leq i\leq
  m\}\] and \[X^T=\{(w,w)\cdot x_i\mid w\in W, 1\leq i\leq m\}.\]

With regards to the existence of filtrable decompositions we show in
the following proposition that the $T\times T$-variety $\mathbb{X}$
admits a filtrable Bialynicki-Birula decomposition if and only if the
toric $T$-variety $X$ also admits one.

\bpropo\label{equivalentfiltrationcondition}
\begin{enumerate}
\item Let $\nu_1,\nu_2\in X_*(T)$ such that the $T\times T$-variety
  $\mathbb{X}$ is filtrable for $\nu=(\nu_1,\nu_2)\in X_*(T\times
  T)$. Then the toric $T$-variety $X$ is filtrable for
  $\nu_1- \nu_2\in X_*(T)$.
\item If the toric $T$-variety $X$ is filtrable for a dominant generic
  $\nu_0 \in X_*(T)$ then the $T\times T$-variety $\mathbb{X}$ is
  filtrable for $\nu=(N\nu_0,\nu_2)$ with $\nu_2$ regular
  anti-dominant and \begin{equation}\tag{\eqref{eq:cond_generic}} N>
    max\left\{\frac{|\langle \mu_{i,\sigma},w\cdot \nu_2\rangle |}{|\langle
        \mu_{i,\sigma},\nu_0\rangle |} \, : \, \sigma \in
      \mathcal{F}(l),\, 1\le i \le q_\sigma \right\}
    \end{equation}
\end{enumerate}
\epropo
\begin{proof}
  {\it Proof of (1):} If
  $\mathbb{X}=Z_1 \supseteq \cdots \supseteq Z_{h}\supseteq
  Z_{h+1}=\emptyset$ is a filtration by closed subsets such that
    \[\forall ~1\le i \le h,\, Z_i\setminus Z_{i+1}= \mathbb{X}_i\] (a
    Bialynicki-Birula cell of $\mathbb{X}$) where $h=m\cdot |W|^2$.
    Let $Z'_i=Z_i \cap X$ for all $i$. The $Z'_i$'s are closed in $X$
    and one has :
    \[X=Z'_{1}\supseteq \cdots \supseteq
      Z'_{h}\supseteq Z'_{h+1}=\emptyset\]

    with
    \[\forall ~ 1\le i\le h,\, ~Z'_i \setminus Z'_{i+1}=
      \mathbb{X}_i\cap X= X_i\] a Bialynicki-Birula cell of $X$ or
    $\emptyset$. In particular, $\mathbb{X}_i\cap X\neq \emptyset$ if
    and only if $\mathbb{X}_i$ is the Bialynicki-Birula cell at a
    $T\times T$-fixed point of the form $(w,w)\cdot x_j$ for some $w\in W$ and
    $1\leq j\leq m$. Hence $X$ is filtrable.

  {\it Proof of (2):} Let \[p : \mathbb{X}\to \overline{G_{ad}}\] the
  $G\times G-$projection onto the wonderful compactification of
  $G_{ad}$. As seen before, the one parameter subgroup $\nu$ in the
  statement is generic for the $T\times T$-variety $\mathbb{X}$ and so
  for $\overline{G_{ad}}=p(\mathbb{X})$ also. Since
  $\overline{G_{ad}}$ is a projective variety, there is a filtration :
\[\overline{G_{ad}}=Y_1\supset \cdots\supset Y_K\supset Y_{K+1}=\emptyset\]
with
\[\forall~~ 1\le i \le K,\, Y_i\setminus Y_{i+1}=
  (\overline{G_{ad}})_i\] a Bialynicki-Birula cell (associated to the
one parameter subgroup $\nu$) centered in a certain $T\times T$-fixed
point $(w_1^i,w_2^i)x_0\in \overline{G_{ad}}$, where
$w^i_1,w^i_2\in W$ and $x_0$ is the only $B\times B^-$-fixed point in
$\overline{G_{ad}}$. In particular, \[K=|W|^2,\] which is the number of
$T\times T$-fixed points in $\overline{G_{ad}}$.

For all $1\le i\le K$,
\[p^{-1}(Y_i\setminus Y_{i+1})=p^{-1} (Y_i) \setminus p^{-1} (Y_{i+1})=
  \bigcup_{\sigma\in \mathcal{F}(l)}\mathbb{X}_{\nu,(w^i_1,w^i_2)\cdot x_\sigma}.\]

Let $1\le i\le K$. Let $w_1:=w^i_1,\, w_2=w^i_2$. Let
$n_1,n_2\in N(T)$ such that $n_1T=w_1,\, n_2T=w_2$. As in the proof of
Proposition \ref{equivalentcellularcondition}, (see \eqref{eq:cellaff}),
there is a $T\times T$-equivariant isomorphism :

\[\varphi \, : \, n_1U \times n_2U^-\times X_0\to (n_1,n_2)\cdot \mathbb{X}_0,\, (u,v,z)\mapsto (u,v)\cdot z\]
which induces isomorphisms for the cells:
\[
  U\cap w_1^{-1} U w_1 \times U^-\cap w_2^{-1}U^-w_2 \times
  {X}_{w_1^{-1} \cdot \nu_1-w_2^{-1} \cdot \nu_2, x_\sigma}\simeq
  \mathbb{X}_{(\nu_1,\nu_2),(w_1,w_2)\cdot x_\sigma }\, ,
\]
for all $\sigma \in \mathcal{F}(l)$.

The action of $(w_1,w_1)$ gives for the cells in the toric variety $X$ :
\begin{equation}
  \forall \sigma\in \mathcal{F}(l),\, (w_1,w_1)\cdot {X}_{w_1^{-1} \cdot \nu_1-w_2^{-1} \cdot \nu_2,x_\sigma}= X_{\nu_1-w_1w_2^{-1} \cdot \nu_2,(w_1,w_1)\cdot x_\sigma}\\
  =X_{\nu_0,(w_1,w_1)\cdot x_\sigma}\end{equation}
according to Lemma \ref{lem:celluparam}.

Since the toric $T$-variety $X$ is supposed to be filtrable for
$\nu_0$, there is a filtration of the toric variety which gives by
restriction a filtration
\[\bigcup_\sigma X_{\nu_0,(w_1,w_1)\cdot x_\sigma} = F_1\supset
  \cdots \supset F_q\supset F_{q +1}=\emptyset\] where the $F_i$'s
are closed subsets of $\displaystyle \bigcup_\sigma X_{\nu_0,(w_1,w_1)\cdot x_\sigma}$
which is a locally closed subset of the toric variety $X$, and
with\[ \forall~~ 1\le j \le m,\, F_j\setminus F_{j+1} =
  X_{\nu_0,(w_1,w_1)\cdot x_j}
\]
for $1\leq j\leq m$. In particular, $q=m=|\mathcal{F}^+(l)|$.

So if one sets 
\[
  \forall ~~1\le j \le m,\,
  Z_{i,j}:=\]\[\overline{\varphi\left(U\cap w_1^{-1} U w_1 \times U^-\cap
      w_2^{-1}U^-w_2 \times (w_1^{-1},w_1^{-1})\cdot
      F_j\right)}\bigcup p^{-1}(Y_{i+1})
\]
(the closure is taken in $\mathbb{X}$) one gets closed subsets of $\mathbb{X}$ such that :
\[
p^{-1}( Y_i)=Z_{i,1}\supset \cdots \supset Z_{i,m}\supset Z_{i,m+1}=p^{-1} (Y_{i+1})
\]
and
\[\forall ~~1\le j \le m,\, Z_{i,j}\setminus Z_{i,j+1}= \]
\[
  \varphi \left(U\cap w_1^{-1} U w_1 \times U^-\cap w_2^{-1}U^-w_2 \times
    (w_1^{-1},w_2^{-1})\cdot {X}_{\nu_0,x_{j}}\right)\]
\[=\mathbb{X}_{(\nu_1,\nu_2),(w_1,w_2)\cdot x_{j} }
\]
That is true for all $1\le i \le K$ so one gets a filtration of
$\mathbb{X}$ by all the \[Z_{i,j},\, 1\le i \le K,\, 1\le j \le m \]
with each successive difference one Bialynicki-Birula cell of
$\mathbb{X}$. So $\mathbb{X}$ is filtrable for $\nu$.
\end{proof}

Combining Proposition \ref{equivalentcellularcondition} and
Proposition \ref{equivalentfiltrationcondition}, by Definition
\ref{cellular} we have the following theorem:

\bth\label{equivcellularcriterion} The toroidal embedding $\mathbb{X}$
is $T\times T$-cellular if and only if the the toric variety $X$ is
$T$-cellular.  \eeth

        {\it Remark.} If $\tT \to T$ is a surjective morphism of torus, then $\mathbb{X}$ is a $\tT\times \tT$-variety and indeed $\mathbb{X}$ is $\tT \times \tT$-cellular (respectively $X$ is $\tT$-cellular) if and only if $\mathbb{X}$ is $T\times T$-cellular (respectively if and only if $X$ is $T$-cellular).

\section{$T\times T$-stable curves in $\mathbb{X}$}\label{invariant curves}

Let $\mathbb{X}$ be a toro\"idal embedding which is cellular for some
one-parameter subgroup $\nu=(\nu_1,\nu_2)$ with $\nu_1,-\nu_2$ regular
dominants in $X_*(T)$. We shall prove that $\mathbb{X}$ satisfies
Assumption \ref{as}.

Let $\mathcal{F}\subset X_{*}(T)$ the fan of the toric variety $X$ associated to $\mathbb{X}$. Let $l$ be the rank of $T$.

\bpropo\label{toroidalassumption}
\begin{enumerate}
    \item Let $C$ be a closed irreducible $T\times T-$stable curve in $\mathbb{X}$. Then $C\simeq  \mathbb{P}^1$.

Let $x\in \mathbb{X}^{T\times T}$ be a $T\times T$-fixed point.
    
  \item  The weights by which  $T\times T$ acts on the outgoing $T\times T$-stable curves of $\mathbb{X}$ at $x$ are pairwise linearly independent.
  \item The number of $T\times T$-stable curves outgoing from $x$
    is $\dim \mathbb{X}^+_{\nu,x}$, the dimension of the
    Bialynicki-Birula cell centered in $x$.
\end{enumerate}

\epropo

\begin{proof}
    \begin{enumerate}
    \item Let $x$ be a $T\times T$-fixed point of $C$. One can suppose
      that $x=x_\sigma\in X^+$ for some maximal cone
      $\sigma\in \mathcal{F}^+(l)$. The intersection
      $C\cap \mathbb{X}_0 \neq \emptyset$. One has a
      $T\times T$-equivariant isomorphism of varieties :
        \[\varphi : U\times U^-\times X_0\to \mathbb{X}_0,\, ((u,v),x)\mapsto (u,v)\cdot x.\]
        Let $U_\sigma=Hom_{sg}(S_{\sigma},\mathbb{C})$ denote the
        $T$-stable open affine subvariety corresponding to the cone
        $\sigma$.  The curve
        \[C':=\varphi^{-1}(C\cap \varphi(U\times U^-\times
          U_\sigma)) \] is $T\times T$-stable in the variety
        $U\times U^-\times U_\sigma$. Now, as an affine toric
        $T-$variety, $U_\sigma$ is a closed $T$-subvariety of some
        $T-$invariant affine space $V_\sigma$ with some $T-$weights
        $\chi \in X^*(T) $ which are pairwise non-proportional. So
        $U\times U^-\times U_\sigma$ is a closed
        $T\times T$-subvariety of a $T\times T$-invariant affine space
        $V=U\times U^-\times V_\sigma$, where $T\times T$ acts by
        weights $(\alpha,0)$, respectively $(0,-\alpha)$, with
        positive roots $\alpha\in \Phi^+$, and by $(\chi,-\chi)$, for
        the weights $\chi\in X^*(T)$ by which $T$ acts on $V_\sigma$.
        All these weights are pairwise non-proportional, so for
        dimension reasons, the curve
        \[ \varphi^{-1}(C\cap \varphi(U\times U^-\times U_\sigma))
        \] is included in some $T\times T$-eigenspace $V_\nu$ for just
        one weight $\nu\in X^*(T\times T)$ of the space $V$.
\begin{itemize}
\item If $\nu=(\alpha,0)$, for some positive root $\alpha$, then
  $C'\subset U_\alpha\times 1\times \{x_\sigma\}$.
\item If $\nu=(0,-\alpha)$, for some positive root $\alpha$, then
  $C'\subset 1\times U_{-\alpha}\times \{x_\sigma\}$.
\item If $\nu=(\chi,-\chi)$ for some $0\neq \chi\in X^*(T)$, then
  $C'\subset 1\times 1\times U_\sigma$.

\end{itemize}
Since $C=\overline{C'}$, in the first case,
$C=\overline{\Big(U_\alpha\times 1\Big)\cdot x_\sigma}$, in the second,
$C=\overline{\Big(1\times U_{-\alpha}\Big)\cdot x_\sigma}$, in the third, $C$ is
a $T-$stable curve of the toric $T$-variety
$X=\overline{U_\sigma}$.

For a positive root $\alpha$, let the rank-one reductive group
$G_\alpha=Z_G(\ker \alpha)^\circ$. The map
$G_\alpha \to \mathbb{X},\, g \mapsto (g,1)\cdot x_\sigma$ induces an
isomorphism
$\mathbb{P}^1\simeq G_\alpha/B^-\cap G_\alpha \simeq
\overline{\Big(U_\alpha\times 1\Big)\cdot x_\sigma}$. For the same reason,
$\mathbb{P}^1\simeq G_\alpha/B\cap G_\alpha \simeq \overline{\Big(1\times
  U_{-\alpha}\Big)\cdot x_\sigma}$. Thirdly, all the $T$-invariant
irreducible curves of a toric $T$-variety are $\mathbb{P}^1$s.

\item Let $x=(w_1,w_2)\cdot x_\sigma$ be a $T\times T$-fixed point in
  $\mathbb{X}$ with $w_1,w_2\in W$ and $\sigma\in \mathcal{F}^+(l)$ a
  maximal cone.

        Let $n_1,n_2\in N(T)$ such that $n_1T=w_1,n_2T=w_2 \in W$.

        The $T\times T$-equivariant open embedding :
\[
U\times U^-\times U_\sigma \to \mathbb{X},\, ((u,v),z)\mapsto (u,v)\cdot z
\]
  induces another $T\times T$-equivariant open embedding :
\[\varphi_{n_1,n_2}:
  n_1U\times n_2U^-\times U_\sigma \to \mathbb{X}\] which maps
  $((n_1u,n_2v),z)\mapsto (n_1u,n_2v)\cdot z .$ Let
$\Omega\subset \mathbb{X}$ be the image of $\varphi_{n_1,n_2}$.

A $T\times T$-invariant outgoing curve $C$ in $\mathbb{X}$ at $x$
is a $T\times T$-stable curve $C$ such that the cell
$C^+_{\nu,x}= C\cap \mathbb{X}^+_{\nu,x}$ is open in $C$, which
means that $C=\overline{C^+_{\nu,x}}$.  So there is a bijection
between the $T\times T$-stable curves in the cell
$\mathbb{X}^+_{\nu,x}$ and the $T\times T$-stable curves in
$\mathbb{X}$ outgoing at $x$.

Now, by the above embedding $\varphi_{n_1,n_2}$, the cell
$\mathbb{X}^+_{\nu,x}$ is isomorphic to :
\[n_1(w_1^{-1} Uw_1\cap U)\times n_2(w_2 ^{-1} U^-w_2\cap U^-)\times
  U_{\sigma,x}^+ \] where $U_{\sigma,x}^+$ is a cell of $U_\sigma$ for
some generic one-parameter subgroup of $T$.

Since we suppose that $\mathbb{X}$ is cellular, the cell
$\mathbb{X}^+_{\nu,x}$ is smooth and so the cell $U_{\sigma,x}^+$ is
an affine $T\times T-$invariant affine space.  The torus $T\times T$
acts
\begin{itemize}
\item on $\varphi_{n_1,n_2}(n_1 U\times n_2\times x_\sigma)$ with
  weights of the form \[(w_1^{-1}(\alpha),0),\, \alpha\in \Phi^+,\]

\item on $\varphi_{n_1,n_2}(n_1\times n_2 U^-\times x_\sigma)$ with
  weights of the form \[(0,w_2^{-1}(\alpha)),\, \alpha\in \Phi^-,\]
\item on $\varphi_{n_1,n_2}(n_1\times n_2\times U_\sigma)$ with
  weights of the form \[(w_1^{-1} (\chi),w_2^{-1}(\chi)),\, \] where
  the weights $\chi\in X^*(T)$ can be chosen among a set of pairwise
  non collinear generators of the semigroup $\sigma^\vee\cap X^*(T)$.
\end{itemize}

Therefore the weights by which $T\times T$ acts on the affine variety 
\[\Omega=\varphi_{n_,n_2}(n_1 U\times n_2 U^-\times U_\sigma) \]
are pairwise non colinear a fortiori those by which $T\times T$ acts
on the cell
\[\mathbb{X}^+_{\nu,x} = \varphi_{n_1,n_2}(n_1(w_1^{-1} Uw_1\cap
  U)\times n_2(w_2 ^{-1} U^-w_2\cap U^-)\times U_{\sigma,x}^+ ) .\]

\item As seen before, the outgoing $T\times T$-stable curves at $x$
  in $\mathbb{X}$ are in bijection with the $T\times T$-stable
  curves in the cell $\mathbb{X}^+_{\nu,x}$. But the cell
  $\mathbb{X}^+_{\nu,x}$ is a $T\times T$-affine space (because
  $\mathbb{X}$ is supposed to be cellular) with $T\times T$-weights
  pairwise non-colinear and with eigenspaces of dimension one.

  The $T\times T$-stable curves of the cell $\mathbb{X}^+_{\nu,x}$
  are precisely the $\varphi_{n_1,n_2}(d)$ where $d$ are the
  eigenlines of the $T\times T$-affine space
\[n_1(w_1^{-1} Uw_1\cap U)\times n_2(w_2 ^{-1} U^-w_2\cap U^-)\times U_{\sigma,x}^+.\]

The number of these eigenlines is in fact  

\[\dim \left(n_1(w_1^{-1} Uw_1\cap U)\times n_2(w_2 ^{-1} U^-w_2\cap
    U^-)\times U_{\sigma,x}^+ \right)=\dim \mathbb{X}^+_{\nu,x} .\]

\end{enumerate}
\end{proof}

In the proposition below we describe the $T\times T$-points  which are linked by an irreducible $T\times T$-stable curve in $\mathbb{X}$.

If $C$ is a $T\times T$-stable curve, if $w_1,w_2\in W$, we note
$(w_1,w_2)C=(n_1,n_2)\cdot C$ for any $n_1,n_2\in N(T)$ such that
$n_iT=w_i,\, i=1,2$. We say that the $T\times T$-stable curves $C$ and
$C'=(w_1,w_2)C$ are $W\times W-$conjugate.

\bpropo\label{descriptioninvariantcurves}
Let $C$ be a closed irreducible $T\times T$-stable curve in $\mathbb{X}$. Then $C$ joins two $T\times T$-fixed points and $C$ is $W\times W-$conjugate to one of the following curves.

\begin{enumerate}
    \item An irreducible $T\times T$-stable curve included in a closed $G\times G-$orbit of $\mathbb{X}$ joining two $T\times T$-fixed points :
      \[x_\sigma,\, (s_\alpha,1)\cdot x_\sigma\] and on which $T\times T$
      acts with a weight $(\pm \alpha,0)$ for some positive root
      $\alpha$ or an irreducible $T\times T$-stable curve included in
      a closed $G\times G-$orbit of $\mathbb{X}$ joining two
      $T\times T-$fixed points :
      \[x_\sigma,\, (1,s_\alpha)\cdot x_\sigma\] and on which
      $T\times T$ acts with a weight $(0,\pm \alpha)$ for some
      positive root $\alpha$.

    \item An irreducible $T\times T-$stable curve not included in any closed $G\times G-$orbit of $\mathbb{X}$, joining two $T\times T$-fixed points 
      \[x_\sigma,\, (s_\alpha,s_\alpha)\cdot x_\sigma\] on which,
      $T \times T$ acts with a weight $\pm(\alpha,-\alpha)$, for some
      $\sigma\in\mathcal{F}^+(l)$ a maximal cone and $\alpha$ a simple
      root orthogonal to a facet of $\sigma$.

    \item An irreducible $T\times T-$stable curve not included in any
      closed $G\times G-$orbit of $\mathbb{X}$ joining two
      $T\times T-$fixed points :
    \[x_\sigma,\, x_{\sigma'}\]
    
    for some $\sigma,\sigma'\in\mathcal{F}^+(l)$ two maximal cones
    which share a facet in the fan $\mathcal{F}^+$ and such that
    $T\times T$ act on $C$ with a weight $(\chi,-\chi)$ with
    $\chi\in X^*(T)$ orthogonal to the facet $\sigma\cap \sigma'$.
\end{enumerate}
\epropo

\begin{proof}
  Let $x\in C^{T\times T}$. There is a $(w_1,w_2)\in W\times W$ such
  that $x=(w_1,w_2)\cdot x_\sigma$ for some maximal cone
  $\sigma\in \mathcal{F}^+$. So $(w_1,w_2)^{-1}C$ is a
  $T\times T-$stable curve which contains $x_\sigma$. So we can
  suppose that $w_1=w_2=1$.

    Then according to the proof of the previous proposition

    \begin{itemize}
    \item Either $C=\overline{(U_\alpha\times 1)\cdot x_\sigma}$ or
      $C=\overline{(1\times U_{-\alpha})\cdot x_\sigma}$ for some positive
      root $\alpha$. In that case the $T\times T$-fixed points of $C$
      are $x_\sigma$ and $(s_\alpha,1)\cdot x_\sigma$ or $x_\sigma$
      and $(1,s_\alpha)\cdot x_\sigma$. This is the case $(1)$ of the
      statement.

    \item Either $C$ is a $T$-stable curve of $\overline{U_\sigma}$
      where $U_\sigma=Hom_{sg}(S_{\sigma},\mathbb{C})$ denotes the
      $T$-stable open affine subvariety of the toric $T$-variety $X$
      associated to the cone $\sigma$.

      In that case, $C=\overline{C\cap U_\sigma}$ is a $T$-stable
      irreducible curve of the toric $T$-variety $\overline{U_\sigma}$
      and so a curve with two $T\times T$-fixed points : $x_\sigma$
      and $x_{\sigma'}$ where $\sigma'$ is a maximal cone in the fan
      $\mathcal{F}$ which share a facet with $\sigma$. If
      $\sigma'\subset \mathcal{F}^+$, we are in the case $(3)$ of the
      statement.

      If $\sigma'\not\subset \mathcal{F}^+$, then there is a
      $1\neq w \in W$ such that $\sigma'\subset w\mathcal{F}^+$. So
      $\sigma\cap \sigma'\subset \mathcal{C}^+\cap w\mathcal{C}^+$. As
      $\sigma\cap \sigma'$ is a cone of dimension $l-1$ in
      $\mathcal{F}$,
      $\dim \mathcal{C}^+\cap w\mathcal{C}^+\ge l-1\Rightarrow
      w=s_\alpha$ for some simple root $\alpha$. So
      $\sigma'=s_\alpha \sigma_1$ for some maximal cone
      $\sigma_1\in \mathcal{F}^+$.  For dimension reasons,
      $\sigma\cap \sigma'= \sigma\cap
      s_\alpha\sigma_1=\alpha^\perp\cap \mathcal{C}^+\Rightarrow
      \alpha^\perp\cap \mathcal{C}^+ \subset \sigma$ (because
      $\alpha^\perp\cap \mathcal{C}^+$ is a cone in
      $\mathcal{F}$). Let $v\in \sigma\setminus \alpha^\perp$. Then
      $\sigma_1\subset {\mathbb{R}}_{\ge 0}v + \alpha^\perp\cap \mathcal{C}^+
      \subset \sigma \Rightarrow \sigma_1\subset \sigma \Rightarrow
      \sigma_1=\sigma$. So $\sigma'=s_\alpha \sigma$ and
      $x_{\sigma'}=(s_\alpha,s_\alpha)\cdot x_\sigma$ with $\alpha $ a
      simple root such that
      $\alpha^\perp\cap \mathcal{C}^+=\sigma\cap \sigma'$. That is the
      case $(2)$.

\end{itemize}
\end{proof}

\brem \label{rem:onecurve} In particular, in the above proof, we have also shown that for
two fixed points $x,x'\in \mathbb{X}^{T\times T}$ there is at most one
$T\times T$-stable irreducible curve in $\mathbb{X}$ which joins $x$
and $x'$.  Indeed, suppose that $x=x_\sigma$ with $\sigma$ a maximal cone in $\mathcal{F}^+(l)$. If $x'\in \mathbb{X}^{T\times T}$, then, if it exists, the curve that joins $x$ and $x'$ is

\begin{itemize}
\item  $\overline{(U_\alpha\times 1)\cdot  x_\sigma}$ if $x'=(s_\alpha,1).x$ for some positive root $\alpha$ ;

\item  $\overline{(1\times U_{-\alpha})\cdot x_\sigma}$ if $x'=(1,s_\alpha).x$ for some positive root $\alpha$ ;
  \item $\overline{V(\tau)}$, closure in the toric variety $X$, if $x'=x_{\sigma'}$ for some maximal cone $\sigma'\in \mathcal{F}(l)$ such that $\sigma \cap \sigma'=\tau$ is a facet of $\sigma$.
\end{itemize}
\erem

\section{Description of $K^0_{T_{comp}}(X)$ and
  $K^0_{T_{comp}}(X^+)$}\label{eqkttvposweyl}

In this section we assume that $X=X(\mathcal{F})$ is $T$-cellular.

We now show that the filtrable Bialynicki-Birula decomposition of
$X$ induces a filtrable Bialynicki-Birula decomposition for $X^+$ as follows:

Recall that the $T$-fixed points of $X$ are $w\cdot x_{\sigma_i}$
correponding to the maximal dimensional cones $ \sigma_i\in \mathcal{F}^+$ for
$1\leq i\leq m$ and $w\in W$.

Let $Y_{w,i}$ denote the Bialynicki Birula cells in $X$ corresponding
to the $T$-fixed point $w\cdot x_{\sigma_i}$ and the one-parameter
subgroup $\lambda_{v}$. Let $\tau_i$ be the face of $\sigma_i$ as
defined in Section \ref{cellulartv}. Then it can be seen that
$Y_i=Y_{1,i}\subseteq X^{+}$  is indeed the Bialynicki Birula cell
in $X^{+}$ corresponding to the $T$-fixed point $x_{\sigma_i}$ and the
one-parameter subgroup $\lambda_{v}$. In particular, if $\sigma_i$ is
a smooth maximal dimensional cone in $\mbox{star}(\tau_i)$ then
$\sigma_i$ is a smooth maximal dimensional cone in
$\mbox{star}(\tau_i)_+$ for $1\leq i\leq m$. Here
$\mbox{star}(\tau_i)_+$ (resp. $\mbox{star}(\tau_i)$) denotes the fan in
$N/N(\tau_i)$ whose maximal dimensional cones correspond to the
maximal dimensional cones in $\mathcal{F}^+$ (resp. $\mathcal{F}$)
containing $\tau_i$.

Thus $Y_i$ for $1\leq i\leq m$ are the Bialynicki-Birula cells of
$X^{+}$ corresponding to its $T$-fixed points $x_{\sigma_i}$ for
$1\leq i\leq m$ and the one-parameter subgroup $\lambda_{v}$. Moreover
since $\displaystyle X=\bigsqcup_{w\in W,1\leq i\leq m}Y_{w,i}$ and
$\displaystyle X^{+}=\bigsqcup_{i=1}^m Y_i$ it follows that
$Y_{w,i}\cap X^{+}=Y_i$ if and only if $w=1$ and
$Y_{w,i}\cap X^{+}=\emptyset$ otherwise.

From the conditions \eqref{star} and \eqref{clstar} it also follows
that the ordering of $X^{T}$ giving a filtrable Bialynicki cellular
decomposition induces an ordering of $(X^{+})^{T}$. Thus without loss of
generality we can assume that the ordering of $X^{T}$ is given by
$\{\sigma_1,\ldots, \sigma_m, \cdots \}$ where the ordering of
$(X^{+})^{T}$ is $\{\sigma_1,\ldots, \sigma_m\}$.  In other words we have
a decreasing filtration:
\[X=Z'_1\supseteq Z'_2\supseteq \cdots \supseteq Z'_m\supseteq
  Z'_{m+1}\cdots \supseteq Z'_{m\cdot |W|}\] and
\[X^{+}=Z_1\supseteq Z_2\supseteq\cdots\supseteq Z_m\supseteq
  Z_{m+1}=\emptyset \] where
$Z_i'\setminus Z_{i+1}'=Z_i\setminus Z_{i+1}=Y_i\simeq
\mathbb{C}^{n-k_i}$ and $Z'_i\cap X^{+}=Z_i$ for $1\leq i\leq m$ and
$Z_i'\cap X^{+}=\emptyset $ for $i\geq m+1$. Here $\dim(\tau_i)=k_i$ for
$1\leq i\leq m$. Since the closure $\overline{Y_i}$ of $Y_i$ in $X$ is
contained in $\displaystyle \bigsqcup_{j\geq i} Y_j$, it implies that
the closure $\overline{Y_i}^+$ of $Y_i$ in $X^{+}$ is contained in
$\displaystyle \bigsqcup_{j\ge i}^m Y_j$. Hence $X$ is $T$-cellular
implies that $X^{+}$ is $T$-cellular.

Henceforth we assume that $X=X(\mathcal{F})$ and hence
$X^{+}=X(\mathcal{F}^+)$ are $T$-cellular.

Thus by Theorem \ref{tveqktdes1} we have the following description of
$K^0_{T_{comp}}(X(\mathcal{F}))$.  \bth\label{main3} The ring
$K^0_{T_{comp}}(X(\mathcal{F}))$ consists in all families
$(a_{\sigma, w})_{\sigma\in \mathcal{F}(l),w\in W}$ of elements in
$R(T_{comp})^{|W|\cdot m}$ such that:

  \begin{enumerate}
  \item[(i)]$a_{\sigma,ws_{\alpha}}\equiv a_{\sigma,w} \pmod{1-e^{-w(\alpha)}}$
    whenever $\alpha\in \Delta$ and $\sigma\in \mathcal{F}^+(l)$ has a
    facet orthogonal to $\alpha$.

  \item[(ii)] $a_{\sigma,w}\equiv a_{\sigma',w} \pmod{1-e^{-w(\chi)}}$
    whenever $\chi\in X^*(T)$ and the cones
    $\sigma,\sigma'\in \mathcal{F}^+(l)$ have a common facet
    orthogonal to $\chi$.
\end{enumerate}    
\eeth Consider the subring $\mathcal{A}^+$ of
$R(T_{comp})^{|\mathcal{F}^+(l)|}$ defined as the
$(a_{\sigma})\in R(T)^{|\mathcal{F}^{+}(l)|}$ such that
$a_{\sigma}\equiv a_{\sigma'} \pmod{(1-e^{-\chi})}$ whenever
$\chi\in X^*(T)$ and the cones $\sigma,\sigma'\in \mathcal{F}^+(l)$
have a common facet orthogonal to $\chi$.  We have the following
theorem (see \cite[Remark 2.4]{u1} for the analogous statement in
algebraic equivariant $K$-theory): \bth\label{cs+} Consider the
restriction map
  \[\iota^*: K^0_{T_{comp}}(X^{+})\lra K^0_{T_{comp}}((X^{+})^{T}).\] Then the image
  $\iota^*(K^0_{T_{comp}}(X^{+}))$ can be identified with the
  $R(T_{comp})$-subalgebra $\mathcal{A}^+$ of
  $K^0_{T_{comp}}((X^{+})^{T})$. \eeth
    \begin{proof} Let
      $(a_{\sigma})_{\sigma\in \cf^{+}(l)}\in \mathcal{A}^+$. Let
      $a_{\sigma,w}:=w\cdot a_{\sigma}$ for each $w\in W$. Then we claim that
      $(a_{\sigma,w})\in R(T_{comp})^{|\mathcal{F}^+(l)|\cdot |W|}$ 
      satisfies (i) and (ii) of Theorem \ref{main3} above. This can be seen as follows: let $\alpha\in \Delta,\, \sigma\in \mathcal{F}^+(l)$ with a facet orthogonal to $\alpha$,
      \[a_{\sigma,w}-a_{\sigma, ws_{\alpha}}=w\cdot a_{\sigma}-w\cdot
        s_{\alpha}\cdot a_{\sigma}.\] In $R(T_{comp})$ it is known
      that
      $a_{\sigma}\equiv s_{\alpha}\cdot a_{\sigma}
      \pmod{(1-e^{-\alpha})}$. Now, translating by $w$ this further
      implies that
      $w\cdot a_{\sigma}\equiv w\cdot s_{\alpha}\cdot a_{\sigma}
      \pmod{(1-e^{-w(\alpha)})}$ for $w\in W$. Moreover, by assumption
      $a_{\sigma}\equiv a_{\sigma'} \pmod{(1-e^{-\chi})}$ whenever
      $\chi\in X^*(T)$ and the cones
      $\sigma,\sigma'\in \mathcal{F}^+(l)$ have a common facet
      orthogonal to $\chi$. Thus for every $w\in W$,
      $w\cdot a_{\sigma}\equiv w\cdot a_{\sigma'}
      \pmod{(1-e^{-w(\chi)})}$ whenever $\chi\in X^*(T)$ and the cones
      $\sigma,\sigma'\in \mathcal{F}^+(l)$ have a common facet
      orthogonal to $\chi$.  By Theorem \ref{main3} this implies that
      $(a_{\sigma,w})_{w\in W, \sigma\in \cf^{+}(l)}$ is in
      $K^0_{T_{comp}}(X)$.

  Now, consider the chain of inclusions
  $(X^{+})^{T}\stackrel{\iota_1}{\hra} X^T\stackrel{\iota_2}{\hra}
  X$. This induces the following chain of maps of equivariant
  $K$-rings
  \[K^0_{T_{comp}}(X)\stackrel{\iota_2^*}{\ra}
    K^0_{T_{comp}}(X^{T})\stackrel{\iota_1^*}{\ra}K^0_{T_{comp}}((X^{+})^{T}).\]
  By the above arguments it can be seen that if
  $(a_{\sigma})_{\sigma\in \cf^{+}(l)}\in \mathcal{A}^+$, then there exists
  $f\in K_{T_{comp}}(X)$ such that
  $\iota_2^*(f)=(a_{\sigma,w})_{w\in W, \sigma\in \cf^{+}(l)}$ and hence
  $\iota_1^*\circ \iota_2^*(f)=(a_{\sigma})_{\sigma\in \cf^{+}(l)}$. Now
  we have the chain of inclusions
  $(X^{+})^{T} \stackrel{\iota}{\hra}X^{+}\stackrel{\iota_3}{\hra} X$
  which induces the following chain of maps of equivariant $K$-rings
  $K^0_{T_{comp}}(X)\stackrel{\iota_3^*}\ra
  K^0_{T_{comp}}(X^{+})\stackrel{\iota^*}{\ra}
  K^0_{T_{comp}}((X^{+})^{T})$. Since
  $\iota_3\circ \iota=\iota_2\circ \iota_1$ it follows that
  $\mathcal{A}^+$ is contained in the image of the map
  $\iota^*\circ \iota_3^*$ and is therefore also contained in the
  image of the map $\iota^*$. By Lemma \ref{localgkm} it is known that
  the image of $\iota^*$ is contained in $\mathcal{A}^+$. Thus the image
  of $\iota^*$ is equal to $\mathcal{A}^+$. Hence the theorem.
    \end{proof}

    Let $PLP(\mathcal{F}^+)$ denote the ring of piecewise Laurent
    polynomial functions on the fan $\cf^+$ (see \eqref{plp}). Since
    $\mathcal{F}^+$ is the fan corresponding to a subdivision of a
    smooth $l$-dimensional cone in $N$, it satisfies \cite[Assumption
    5.5]{u3} (see \cite[Section 5.2, Remark 7.2]{u3}). Indeed,
    consider $S_{N}:=N_{\mathbb{R}}\setminus \{0\}/\mathbb{R}_{>0}$
    which is homeomorphic to the $(l-1)$-sphere and the projection
    $\uppi: N_{\mathbb{R}}\setminus \{0\}\ra S_{N}$. The polytopal
    complex formed by $\uppi(\sigma\setminus\{0\})$ for
    $\sigma\in \mathcal{F}^+$ is a polytopal subdivision of an
    $(l-1)$-ball and is therefore pure and strongly connected. Thus by
    Theorem \ref{cs+} we get the following result: \bth\label{main2}
    The ring $K^0_{T_{comp}}(X^+)$ is isomorphic to
    $PLP(\mathcal{F}^+)$ as an $R(T_{comp})$-algebra.  \eeth We shall
    skip the proof of the above theorem which follows exactly along
    the same lines as that for a complete $T$-cellular toric variety
    $X=X(\Sigma)$ in \cite[Theorem 5.6]{u3} (see Theorem
    \ref{tveqkthdes2} above), with $\mathcal{A}$ replaced by
    $\mathcal{A}^+$ and $PLP(\Sigma)$ replaced by $PLP(\cf^+)$.

\brem\label{exttorus} In Theorem \ref{main3}, Theorem \ref{cs+} and Theorem \ref{main2} we can replace $T_{comp}$ by $\tT_{comp}$ by considering the action of $\tT$ (respectively $\tT_{comp}$) on $X$ and $X^+$ by means of its surjection to $T$ (respectively $T_{comp}$). 
\erem

\section{Equivariant topological $K$-ring of cellular toroidal
  embeddings}\label{eqtopcelltoroidal}

We shall follow the notations of Section \ref{Introduction} and Section \ref{celltoremb}.

Recall from \cite[Corollary 3.7]{I} that there exists an exact
sequence \be\label{exactseq}1\ra Z\ra \tG:=\tC\times
G^{ss}\stackrel{\uppi}{\ra} G\ra 1\ee where $Z$ is a finite central
subgroup, $\tC$ is a torus and $G^{ss}$ is semisimple and simply
connected. In particular, $\tG$ is {\it factorial} and
$\tB:=\uppi^{-1}(B)$ and $\tT:=\uppi^{-1}(T)$ are respectively a Borel
subgroup and a maximal torus of $\tG$.

Let $\tG_{comp}$ denote a maximal compact subgroup of $\tG$ such that
$\tT_{comp}=\tT\cap \tG_{comp}$ is the maximal compact torus of
$\tT$. Since $\tG$ is factorial, $\pi_1(\tG_{comp})$ is torsion free and hence free.

We shall consider the actions of ${\tG_{comp}}\times \tG_{comp}$ and
${\tT}_{comp}\times \tT_{comp}$ on $\mathbb{X}$ via the projections to $G_{comp}\times G_{comp}$
and $T_{comp}\times T_{comp}$ respectively.

For $\sigma\in \mathcal{F}^+(l)$ we denote by $\mathcal{Z}_{\sigma}\cong G/B^{-}\times G/B$ the corresponding closed orbit with base point  $z_{\sigma}$ and by

$\iota_{\sigma}: K^0_{\tT_{comp}\times \tT_{comp}}(\mathbb{X})\ra
K^0_{\tT_{comp}\times \tT_{comp}}(\cz_{\sigma})$ the corresponding
restriction map. Moreover, for
$f\in K_{\tT_{comp}\times \tT_{comp}}(\cz_{\sigma})$ we denote by
$f_{w_1,w_2}$ the restriction to the point $(w_1,w_2)z_{\sigma}$

Note here that $R(\tT_{comp}\times \tT_{comp})=\mathbb{Z}[X^*(\tT_{comp}\times \tT_{comp})]$ where $e^{(\chi_1,\chi_2)}$ corresponds to the one-dimensional representation where $\tT_{comp}\times \tT_{comp}$ acts via the character $(\chi_1,\chi_2)$ where $(\chi_1,\chi_2)(t_1,t_2)=\chi_1(t_1)\cdot \chi_2(t_2)$.

With the above notations,  the following theorem describes the
$\tT_{comp}\times \tT_{comp}$-equivariant topological $K$-ring of a
cellular toroidal embedding. This generalizes the description of the
$T\times T$-equivariant Chow ring in \cite{Br2} and the
$T\times T$-equivariant Grothendieck ring of a regular
compactification of $G$ to the setting of topological equivariant
$K$-theory.

\bth\label{torus-equivariant K-ring}
For a $T\times T$-cellular toroidal embedding $\mathbb{X}$ of $G$ the map 

\[\prod_{\sigma\in \mathcal{F}_{+}(l)}\iota_{\sigma}:
  K^0_{\tT_{comp}\times \tT_{comp}}(\mathbb{X})\ra \prod_{\sigma\in
    \mathcal{F}_{+}(l)} K^0_{\tT_{comp}\times \tT_{comp}}
  (\cz_{\sigma})\]

is injective and its image consists of  all families $(f_{\sigma})_{\sigma\in \mathcal{F}_{+}(l)}$ of elements in $\prod_{\sigma\in \mathcal{F}_{+}(l)}K^0_{\tT_{comp}\times \tT_{comp}} (\cz_{\sigma})$ such that:

(i)
$f_{\sigma,w_1\cdot s_{\alpha}, w_2\cdot s_{\alpha}}\equiv f_{\sigma,
  w_1, w_2}\pmod {(1-e^{(-w_1(\alpha), w_2(\alpha))})}$ whenever
$\alpha\in \Delta$ and the cone $\sigma\in \mathcal{F}^+(l)$ has a
facet orthogonal to $\alpha$ and

(ii) $f_{\sigma,w_1,w_2}\equiv f_{\sigma',w_1,w_2} \pmod{(1-e^{(-w_1(\chi),w_2(\chi))}})$ whenever the cones $\sigma$ and 
$\sigma'$ have a common facet orthogonal to $\chi\in X^*(T)$. 

\eeth

\begin{proof}

  Since $\mathbb{X}$ is $T\times T$-cellular it can be seen that it is $\tT\times \tT$-cellular for the induced action of $\tT\times \tT$. Hence it
  satisfies Assumption \ref{as} by Proposition
  \ref{toroidalassumption}. Furthermore, we have a precise description
  of the $T\times T$-stable and hence
  $\tT\times \tT$-stable curves in Proposition
  \ref{descriptioninvariantcurves}. By Proposition \ref{localization}
  and Theorem \ref{main1} and
  $K^0_{\tT_{comp}\times \tT_{comp}}(\mathbb{X})$ is isomorphic as a
  $R(\tT_{comp}\times \tT_{comp})$-subalgebra of
  \[\prod_{\sigma\in \mathcal{F}_{+}(l), w_1,w_2\in W} K_{\tT_{comp}\times
      \tT_{comp}}((w_1,w_2)\cdot z_{\sigma})\] given by the relations
  corresponding to the curves listed in Proposition
  \ref{descriptioninvariantcurves}. The first type of curves lie in
  the closed orbit $\cz_{\sigma}\cong G/B^-\times G/B$ and hence
  define its image in
  $\displaystyle \prod_{\sigma\in \mathcal{F}_{+}(l), w_1,w_2\in W}
  K^0_{\tT_{comp}\times \tT_{comp}}((w_1,w_2)\cdot z_{\sigma})$ (see for
  example \cite[Theorem 1.6]{ml}). The second and the third types of
  curves give the following relations respectively:

(i) $f_{\sigma, w_1\cdot s_{\alpha}, w_2\cdot s_{\alpha}}\equiv f_{\sigma, w_1,w_2}\pmod{(1-e^{(-w_1(\alpha),w_2(\alpha))}) }$

(ii) $f_{\sigma, w_1,w_2}\equiv f_{\sigma',w_1,w_2}\pmod {(1-e^{(-w_1(\chi), w_2(\chi))})}$

Since the restriction 
\[\iota :K^0_{\tT_{comp}\times \tT_{comp}}(\mathbb{X})\ra
  K^0_{\tT_{comp}\times \tT_{comp}}(\mathbb{X}^{T\times T})\] is
injective and $\iota$ factors through
$\displaystyle\prod_{\sigma\in \mathcal{F}^+(l)} K^0_{\tT_{comp}\times
  \tT_{comp}} (\cz_{\sigma})$ it follows that
$\displaystyle \prod_{\sigma\in \mathcal{F}^+(l)} \iota_{\sigma}$ is
injective.  Furthermore, the relations (i) and (ii) give the
congruences satisfied by the image of
$\displaystyle\prod_{\sigma\in \mathcal{F}_{+}(l)}\iota_{\sigma}$.
Hence the theorem.

\end{proof}

Since $\mathbb{X}$ and $G/B^-\times G/B$ are both $T\times T$-cellular
$G\times G$-varieties we shall apply Proposition \ref{W-action
  general} below.
\[\Big(K_{\tT_{comp}\times \tT_{comp}}^0(G/B^-\times G/B)\Big)^{W\times
    W}=K^0_{\tG_{comp}\times \tG_{comp}}(G/B^-\times G/B)\] and
\[\Big(K_{\tT_{comp}\times \tT_{comp}}^0(\mathbb{X})\Big)^{W\times
    W}=K^0_{\tG_{comp}\times \tG_{comp}}(\mathbb{X})\] where
$W\times W$ acts on
$K_{\tT_{comp}\times \tT_{comp}}^0(G/B^-\times G/B)$ as follows:

\[((w_1,w_2)\cdot f)\mid_{(w_1', w_2')}=(w_1,w_2)\cdot [f\mid_{(w_1^{-1}\cdot w_1',
  w_2^{-1}\cdot w_2')}]\] and on
$K_{\tT_{comp}\times \tT_{comp}}^0(\mathbb{X})$ as follows:

\[((w_1,w_2)\cdot f)\mid_{\sigma, w_1', w_2'}=(w_{1}, w_{2})\cdot
  [f\mid_{\sigma, w_1^{-1}\cdot w_1', w_2^{-1}\cdot w_2'}]\] for
$(w_1,w_2), (w_1',w_2')\in W\times W$ and $\sigma\in \cf^{+}(l)$.

Furthermore, the maps $\iota$ and
$\displaystyle \prod_{\sigma\in \cf^{+}(l)} \iota_{\sigma}$ are
$W\times W$-equivariant where the $W\times W$ action on
$\displaystyle\prod_{\sigma\in \cf^{+}(l), (w_1, w_2)\in W\times W}
R(\tT_{comp}\times \tT_{comp})$ is given by
\[(w_1,w_2)\cdot (a_{(w,w')})_{(w,w')\in W\times W}=(b_{(w,
    w')})_{(w,w')\in W\times W}\] where
$b_{(w,w')}=(w_1,w_2)\cdot a_{(w_1^{-1}\cdot w, w_2^{-1}\cdot w')}$
for $(w,w')\in W\times W$.

\subsection{\boldmath Description of
  $K^0_{\tG_{comp}\times
    \tG_{comp}}(\mathbb{X})$}\label{descriptioneqtop}

\bcor\label{$G$-equivariant K-ring} The ring
$K^0_{\tG_{comp}\times \tG_{comp}}(\mathbb{X})$ consists in all
families $(f_{\sigma})_{\sigma\in \cf^{+}(l)}$ of elements in
$R(\tT_{comp})\otimes R(\tT_{comp})$ such that

(i) $(s_{\alpha},s_{\alpha})f_{\sigma}\equiv
f_{\sigma}\pmod{1-e^{(-\alpha,\alpha)}}$ whenever $\alpha\in \Delta$
and the cone $\sigma\in \cf^{+}(l)$ has a facet orthogonal to $\alpha$; and

(ii) $f_{\sigma}\equiv f_{\sigma'}\pmod{(1-e^{(-\chi,\chi)})}$ whenever $\chi\in X^*(T)$ and the cones $\sigma$ and $\sigma'$ have a common facet orthogonal to $\chi$.

\ecor
\begin{proof}
  If
  $f\in K_{\tG_{comp}\times
    \tG_{comp}}^0(\mathbb{X})=K_{\tT_{comp}\times
    \tT_{comp}}^0(\mathbb{X})^{W\times W}$ then 

  \[f \mid_{\sigma, w_1', w_2'}=((w_1,w_2)\cdot f) \mid_{\sigma, w_1',
      w_2'}=(w_{1}, w_{2})\cdot [f\mid_{\sigma, w_1^{-1}\cdot w_1',
      w_2^{-1}\cdot w_2'}]\] for $(w_1,w_2), (w_1',w_2')\in W\times W$
  and $\sigma\in \cf^{+}(l)$ (see Proposition \ref{W-action
    general}). Thus
  \[f\mid_{\sigma, w_1^{-1}\cdot w_1', w_2^{-1}\cdot
      w_2'}=(w^{-1}_1,w^{-1}_2) f \mid_{\sigma, w_1', w_2'}\] for
  $(w_1,w_2), (w_1',w_2')\in W\times W$ and $\sigma\in \cf^{+}(l)$.

  Thus the relations (i) and (ii) in Theorem \ref{torus-equivariant K-ring} reduce respectively to

  (i') $(w_1,w_2)(s_{\alpha}, s_{\alpha})f_{\sigma,1,1}\equiv (w_1,
  w_2)f_{\sigma, 1,1}\pmod {1-e^{(-w_1(\alpha),-w_2(-\alpha))}}$

  (ii') $(w_1,w_2)f_{\sigma,1,1}\equiv (w_1, w_2)f_{\sigma', 1,1}\pmod
  {1-e^{(-w_1(\chi),-w_2(-\chi))}}$ for $(w_1,w_2)\in W\times W$.

 % K^0_{\tT_{comp}\times
   % \tT_{comp}}(Z_{\sigma})^{W\times W}
  Since
  \begin{align*} K^0_{\tG_{comp}\times \tG_{comp}}(\cz_{\sigma})=&K^0_{\tG_{comp}\times
    \tG_{comp}}(\tG_{comp}\times \tG_{comp}\times_{\tT_{comp}\times
    \tT_{comp}}z_{\sigma})\\=&K^0_{\tT_{comp}\times
    \tT_{comp}}(z_{\sigma})\cong R(\tT_{comp}\times \tT_{comp}),\end{align*} we
can identify $f_{\sigma}$ with $f_{\sigma, 1,1}$. Thus the only
nontrivial relations are given by:

  (i) $(s_{\alpha},s_{\alpha})f_{\sigma}\equiv
f_{\sigma}\pmod{1-e^{(-\alpha,\alpha)}}$ whenever $\alpha\in \Delta$
and the cone $\sigma\in \cf^{+}(l)$ has a facet orthogonal to $\alpha$; and

(ii) $f_{\sigma}\equiv f_{\sigma'}\pmod{(1-e^{(-\chi,\chi)})}$ whenever
$\chi\in X^*(T)$ and the cones $\sigma$ and $\sigma'$ have a common
facet orthogonal to $\chi$.
\end{proof}

Recall (see \cite{u1,u2}) that we have a split exact sequence:
\be\label{diag} 1\ra \mbox{diag} (\tT)\ra \tT\times \tT\ra \tT\ra 0,\ee where the
second map sends $(t_1,t_2)\ra t_1t_2^{-1}$ and the splitting is given by $t\mapsto (t,1)$.

In particular, $(\alpha, -\alpha)$ and $(\chi, -\chi)$ are characters
on $T\times T$ which are trivial on $\mbox{diag}(T)$ and hence give
rise to a character of $T\times 1$.

\bcor\label{$G$-equivariant K-ring change of variables}
The ring $K_{\tG_{comp}\times
  \tG_{comp}}^0(\mathbb{X})$ consists of all families $(f_{\sigma})_{\sigma \in \mathcal{F}^+(l)}$ of elements in $R(\tT_{comp}\times \{1\})\otimes R(\mbox{diag}(\tT_{comp}))$ such that

(i) $(1,s_{\alpha})f_{\sigma}(u,v)\equiv f_{\sigma}(u,v)\pmod {1-e^{-\alpha(u)}}$ whenever $\alpha\in \Delta$ and the cone $\sigma\in \cf^{+}(l)$ has a facet orthogonal to $\alpha$,

(ii) $f_{\sigma}\equiv f_{\sigma'}\pmod{1-e^{-\chi(u)}}$ whenever
$\chi\in X^*(T)$ and the cones $\sigma$ and $\sigma'$ in $\cf^{+}(l)$
have a common facet orthogonal to $\chi$.

Here the variables $u$ and $v$ correspond to
$R(\tT_{comp}\times \{1\})$ and $R(\mbox{diag}(\tT_{comp}))$
respectively.
\ecor

\begin{proof} It can be seen that the split exact sequence
  \eqref{diag} gives an identification of
  $X^*(\tT_{comp}\times \tT_{comp})$ with
  $X^*(\tT_{comp}\times 1) \oplus X^*(\mbox{diag}(\tT_{comp}))$, where
  the character $(\chi_1,\chi_2)\in X^*(\tT_{comp}\times \tT_{comp})$
  maps to the character $(\chi_1,\chi_2) | _{\tT_{comp}\times 1} \oplus (\chi_1,\chi_2) | _{\mathrm{diag}(\tT_{comp})} =(\chi_1(u), (\chi_1+\chi_2)(v))\in
  X^*(\tT_{comp}\times 1) \oplus X^*(\mbox{diag}(\tT_{comp}))$. Here
  the variables $u$ and $v$ respectively correspond to the factors $X^*(\tT_{comp}\times 1)$ and
  $X^*(\mbox{diag}(\tT_{comp}))$  of
  $X^*(\tT_{comp}\times 1)\oplus X^*(\mbox{diag}(\tT_{comp}))$. In particular, $(-\chi,\chi)$ and $(-\alpha,\alpha)$
  map respectively to $(-\chi(u),0)$ and $(-\alpha(u),0)$.

  We can identify $R(\tT_{comp}\times \tT_{comp})$ with
  \[\mathbb{Z}[e^{\chi_1(u)}\otimes e^{(\chi_1+\chi_2)(v)}:
  (\chi_1,\chi_2)\in X^*(\tT_{comp}\times \tT_{comp})].\] Furthermore,
  note that
  \[(s_{\alpha},s_{\alpha})\cdot
    f_{\sigma}-f_{\sigma}=(s_{\alpha},1)\cdot (1,s_{\alpha})\cdot
    f_{\sigma}-(1,s_{\alpha})f_{\sigma}+
    (1,s_{\alpha})f_{\sigma}-f_{\sigma}\] and
  $(s_{\alpha},1)\cdot f-f\equiv 0\pmod{1-e^{-\alpha(u)}}$ for all
  $f\in R(\tT_{comp}\times \tT_{comp})$. The proof now follows from
  Corollary \ref{$G$-equivariant K-ring}.
\end{proof}

In particular, the image of
$K_{\tG_{comp}\times \tG_{comp}}^0(\mathbb{X})$ under $\iota$ maps
injectively into
$\prod_{\sigma\in \cf^{+}(l)} K^0_{\tT_{comp}\times
  \tT_{comp}}(z_{\sigma})=K^0_{\tT_{comp}\times
  \tT_{comp}}((X^+)^{T})$. Thus $\iota$ restricted to
$K_{\tG_{comp}\times \tG_{comp}}^0(\mathbb{X})$ factors through
$K^0_{\tT_{comp}\times \tT_{comp}}(X^+)$. Since
$\mbox{diag}(\tT_{comp})$ acts trivially on $X^+$. Thus
\[K^0_{\tT_{comp}\times \tT_{comp}}(X^+)\cong
K^0_{\tT_{comp}\times \{1\}}(X^+)\otimes
R(\mbox{diag}(\tT_{comp})).\]

\bpropo\label{chain of inclusions} There is a chain of injective
morphisms of $R(\tG_{comp})\otimes R(\tG_{comp})$-algebras:
\be\label{chain} K^0_{\tT_{comp}}(X^+)\otimes R(\tG_{comp})\subseteq
K_{\tG_{comp}\times \tG_{comp}}^0(\mathbb{X})\subseteq
K^0_{\tT_{comp}}(X^+)\otimes R(\tT_{comp}).\ee Moreover, the
$R(\tT_{comp})$-algebra structure on $K^0_{\tT_{comp}}(X^{+})$ induces
an $R(\tT_{comp})\otimes R(\tG_{comp})$-algebra structure on
$K^0_{\tG_{comp}\times \tG_{comp}}(\mathbb{X})$.  \epropo

\begin{proof}
  Recall from Theorem \ref{cs+} that the canonical restriction
  homomorphism
  \[K^0_{\tT_{comp}}(X^+) \ra K^0_{\tT_{comp}}((X^+)^{T})=
    R(\tT_{comp})^{|\cf^{+}(l)|}\] is injective and the image of
  $K^0_{\tT_{comp}}(X^+)$ in
  $\displaystyle R(\tT_{comp})^{|\cf^{+}(l)|}$ consists of
  $(a_{\sigma})$ such that
  $a_{\sigma}\equiv a_{\sigma'}\pmod{1-e^{-\chi}}$ whenever
  $\chi\in X^*(T)$ and $\sigma$ and $\sigma'\in \cf^{+}(l)$ have a
  common facet orthogonal to $\chi$. Furthermore, the above map is
  compatible with the canonical $R(\tT_{comp})$-algebra structure on
  $K^0_{\tT_{comp}}(X^+)$ and the $R(\tT_{comp})$ algebra structure on
  $\displaystyle R(\tT_{comp})^{|\cf^{+}(l)|}$ given by the diagonal
  map.

  Furthermore, since $R(\tG_{comp})=R(\tT_{comp})^{W}$, it follows from relations
  (i) and (ii) of Corollary \ref{$G$-equivariant K-ring change of
    variables} that $K^0_{\tT_{comp}}(X^+)\otimes R(\tG_{comp})$
  is a subring of $K^0_{\tG_{comp}\times \tG_{comp}}(\mathbb{X})$. 

  In particular, $K^0_{\tG_{comp}\times \tG_{comp}}(\mathbb{X})$ is an
  algebra over $K^0_{\tT_{comp}}(X^+)\otimes R(\tG_{comp})$ and
  hence over $R(\tT_{comp})\otimes R(\tG_{comp})$. Thus we get the
  chain \eqref{chain} of inclusions of
  $R(\tT_{comp})\otimes R(\tG_{comp})$-algebras.  
\end{proof}

Let
\[PL(\cf^{+}):=\left\{ h:\mathcal{C}^+\ra \mathbb{R}\mid \forall \sigma\in
\cf^{+}(l),\, \exists h_\sigma\in X^*(T),\, \forall v\in \sigma,\, h(v)=\langle h_{\sigma},
v\rangle\right\}\] be the set of
piecewise linear functions on $\mathcal{C}^+$ that are linear on the
subdivision $\cf^{+}$.  By the parametrization of line bundles on
spherical varieties (see \cite{Br3}, \cite[Section 1.2.2]{tc}) we know
that the group of isomorphism classes of $\tG\times \tG$-linearized
line bundles on $\mathbb{X}$ is isomorphic to $PL(\cf^{+})$. Let
$\mathcal{L}_{h}$ be the $\tG\times \tG$-linearized line bundle on
$\mathbb{X}$ corresponding to $h\in PL(\cf^{+})$. We consider
$\mathcal{L}_h$ as a $\tG_{comp}\times \tG_{comp}$-linearized line
bundle on $\mathbb{X}$ with the restricted action.

We have a canonical map $PL(\cf^{+})\ra PLP(\cf^+)$ given by
$h\mapsto (e^{h_{\sigma}})_{\sigma\in \cf^+(l)}$ for
$h=(h_{\sigma})_{\sigma\in \cf^+(l)}$ where $PLP(\cf^+)$ denotes the
ring of piecewise Laurent polynomial functions on $\cf^+$ (see \eqref{plp})).

The following theorem is the extension of \cite[Theorem 2.1]{u2} to
the setting of topological equivariant $K$-ring of toroidal
embeddings. We skip the proof as it follows exactly along the same
lines as that of \cite[Theorem 2.1]{u2} using Corollary
\ref{$G$-equivariant K-ring change of variables}, Proposition
\ref{chain of inclusions} and \eqref{dec1}.

\bth\label{decomp} The ring
$K^0_{\tG_{comp}\times \tG_{comp}}(\mathbb{X})$ has the following
direct sum decomposition as a
$K^0_{\tT_{comp}}(X^+)\otimes R(\tG_{comp})$-module: \be\label{dec2}
K^0_{\tG_{comp}\times \tG_{comp}}(\mathbb{X})=\bigoplus
\prod_{\alpha\in I}(1-e^{\alpha(u)})\cdot K^0_{\tT_{comp}}(X^+)\otimes
R(\tT_{comp})_{I}.\ee This direct sum is a free
$K^0_{\tT_{comp}}(X^+)\otimes R(\tG_{comp})$-module of rank
$|W|$ with basis
$\Big\{\prod_{\alpha\in I}(1-e^{\alpha})\otimes f_{v}:v\in C^{I} ~and
~I\subseteq \Delta\Big\}$

where $C^{I}$ and $f_{v}$ are as defined in Section
\ref{free}. Moreover, the component
$K^0_{\tT_{comp}}(X^+)\otimes 1\subseteq
K^0_{\tT_{comp}}(X^+)\otimes R(\tT_{comp})^{W}$ of the direct sum
decomposition can be identified with the subring of
$K^0_{\tG_{comp}\times \tG_{comp}}(\mathbb{X})$ generated by the
classes of ${\tG_{comp}\times \tG_{comp}}$-equivariant line bundles on
$\mathbb{X}$.

\eeth

Since $\displaystyle \prod_{\alpha\in I}(1-e^{\alpha})$ is not a zero divisor in
$R(\tT_{comp})$ and hence in $K_{\tT_{comp}}^0(X^{+})$ which is a free
module over $R(\tT_{comp})$. Thus each piece
$\displaystyle\prod_{\alpha\in I}(1-e^{\alpha(u)})\cdot
K_{\tT_{comp}}^0(X^+)\otimes R(\tT_{comp})_{I}$ is isomorphic to
$K_{\tT_{comp}}^0(X^+)\otimes R(\tT_{comp})_{I}$ via an isomorphism
which maps $\displaystyle\prod_{\alpha\in I}(1-e^{\alpha(u)})\cdot a\otimes b$ to the element
$a\otimes b$.

We have the following theorem which is a generalization of
\cite[Theorem 2.2]{u2} to $\tG_{comp}\times \tG_{comp}$-equivariant
topological $K$-ring of toroidal emebddings. We skip the proof which
involves similar arguments as that of \cite[Theorem 2.2]{u2} to avoid
repetition.

\bth\label{multstr} We have the following decomposition of
$K^0_{\tG_{comp}\times \tG_{comp}}(\mathbb{X})$ as
$K_{\tT_{comp}}^0(X^+)\otimes R(\tG_{comp})$ submodules of
$K_{\tT_{comp}}^0(X^+)\otimes R(\tT_{comp})$
\be\label{dec4}
K^0_{\tG_{comp}\times \tG_{comp}}(\mathbb{X})=\bigoplus
K^0_{\tT_{comp}}(X^+)\otimes R(\tT_{comp})_{I}.\ee

Under the above decomposition, the elements $1\otimes f_v$ for $v\in W$
form a basis for $K^0_{\tG_{comp}\times \tG_{comp}}(\mathbb{X})$ as
$K_{\tT_{comp}}^0(X^+)\otimes R(\tG_{comp})$-module with the
multipication rule given by

\begin{align*} &(1\otimes f_{v})\cdot (1\otimes f_{v'}) \\ & =
  \sum_{J\subseteq I\cup I'}\sum_{w\in C^{J}}\Big(\prod_{\alpha\in
    I\cap I'} (1-e^{\alpha(u)})\cdot \prod_{\alpha\in I\cup
    I'\setminus J} (1-e^{\alpha(u)})\Big) \otimes a_{v,v'}^{w}\cdot
  (1\otimes f_{w})\end{align*} \eeth for $v\in C^{I}$ and
$v'\in C^{I'}$ for $I,I'\subseteq \Delta$ and
$a_{v,v'}^{w}\in R(\tG_{comp})$ are as in \eqref{mult1}.

\subsection{Relation to the \boldmath $\tG_{comp}\times \tG_{comp}$-equivariant
  $K$-ring of the wonderful compactification}.

Recall from \eqref{exactseq} that $\tG=\tC\times G^{ss}$, where $\tC$ is a torus and $G^{ss}$ is semisimple and simply connected. Moreover, $G^{ss}$ is the universal cover of $G_{ad}$.

Recall from Section \ref{toroidal} that we have a canonical surjective morphism $p:\mathbb{X}\ra \overline{G_{ad}}$  which further restricts to a proper surjective morphism $p\mid_{X^+}:X^+\ra \bar{T_{ad}}^+=\mathbb{A}^r$ of toric varieties. Furthermore, $p$ is equivariant with respect to the action of $\tG\times \tG$ where $\tG\times \tG$ acts on $\overline{G_{ad}}$ by means of the quotient $G^{ss}\times G^{ss}$. Thus $p$ is also equivariant under the restricted action of $\tG_{comp}\times \tG_{comp}$ on $\mathbb{X}$ and $\overline{G_{ad}}$.

Thus the pull-back map $p^*: K^0_{\tG_{comp}\times \tG_{comp}}(\overline{G_{ad}})\ra K^0_{\tG_{comp}\times \tG_{comp}}(\mathbb{X})$ induces a canonical $K^0_{\tG_{comp}\times \tG_{comp}}(\overline{G_{ad}})$-structure on $K^0_{\tG_{comp}\times \tG_{comp}}(\mathbb{X})$.

We have the following theorem which describes
$K^0_{\tG_{comp}\times \tG_{comp}}(\mathbb{X})$ as an
$K^0_{\tG_{comp}\times \tG_{comp}}(\bar{G_{ad}})$-algebra.

\bth\label{relwond} We have the following isomorphism of
$R(\tG_{comp})\otimes R(\tG_{comp})$-algebras
\[K^0_{\tG_{comp}\times \tG_{comp}}(\mathbb{X})\cong
  K^0_{\tG_{comp}\times
    \tG_{comp}}(\bar{G_{ad}})\otimes_{R(\tT_{comp})}
  K^0_{\tT_{comp}}(X^+)\] where
$K^0_{\tG_{comp}\times \tG_{comp}}(\bar{G_{ad}})$ is an
$R(\tT_{comp})$-algebra via the map which sends $e^{\chi}$ for
$\chi\in X^*(\tT_{comp})$ to $[\mathcal{L}_{\chi}]$ where
$\mathcal{L}_{\chi}$ is the associated
$\tG_{comp}\times \tG_{comp}$-linearized line bundle on $\bar{G_{ad}}$
(see \cite{dp}). \eeth 

\begin{proof}
 % Since $\bar{G_{ad}}$ is a smooth projective variety we have an
 % isomorphism of the canonical map
 % $\mathcal{K}^0_{\tG\times \tG}(\bar{G_{ad}})\cong
 % K^0_{\tG_{comp}}(\bar{G_{ad}})$ is an isomorphism (see \cite[Lemma
 % 1.12, Remark 1.14]{u1} and Section \ref{eqalgkth}).

  When $\mathbb{X}=\bar{G_{ad}}$, $X^+\cong\mathbb{A}^r$ so that
  \eqref{dec2} gives the decomposition

  \be\label{dec3}  K^0_{\tG_{comp}\times \tG_{comp}}(\bar{G_{ad}})=\bigoplus_{I\subseteq \Delta}
\prod_{\alpha\in I}(1-e^{\alpha(u)})\cdot R(\tT_{comp})\otimes
R(\tT_{comp})_{I}.  \ee
  
  This decomposition is analogous to that in \cite[Theorem 3.3]{u1}, for the
  $\tG\times \tG$-equivariant algebraic $K$-ring of $\bar{G_{ad}}$
  as a $R(\tT)\otimes R(\tG)$-module.

  Combining \eqref{dec2} and \eqref{dec3} we get:
  \[K^0_{\tG_{comp}\times \tG_{comp}}(\mathbb{X})=
    K^0_{\tG_{comp}\times
      \tG_{comp}}(\bar{G_{ad}})\bigotimes_{R(\tT_{comp})\otimes 1}
    K^0_{\tT_{comp}}(X^+)\otimes 1\] where the
  $R(\tT_{comp})\otimes 1$-algebra structure on
  $K^0_{\tG_{comp}\times \tG_{comp}}(\bar{G_{ad}})$ is via the
  decomposition \eqref{dec3}.  It is known by \cite[Theorem 3.6]{u1} that
  $R(\tT)\otimes 1=R(\tT_{comp})\otimes 1$ is the subring of
  $K^0_{\tG_{comp}\times \tG_{comp}}(\bar{G_{ad}})$ generated by
  $Pic_{\tG\times \tG}(\bar{G_{ad}})$ where
  $e^{\chi}\in R(\tT_{comp})$ corresponds to the isomorphism class
  of the line bundle $\mathcal{L}_{\chi}$ on $\bar{G_{ad}}$, which
  has a $\tG\times \tG$-linearization and hence a
  $\tG_{comp}\times \tG_{comp}$-linearization.
\end{proof}

The following corollary describes $K^0_{\tG_{comp}\times \tG_{comp}}(\mathbb{X})$ as an extension of scalars of the $R(\tT_{comp})$-algebra $PLP(\cf^+)$ to the ring $K^0_{\tG_{comp}\times \tG_{comp}}(\overline{G_{ad}})$.

\bcor\label{relwondcor} We have the following isomorphism of
$R(\tG_{comp})\otimes R(\tG_{comp})$-algebras
\[K^0_{\tG_{comp}\times \tG_{comp}}(\mathbb{X})\simeq
  K^0_{\tG_{comp}\times
    \tG_{comp}}(\bar{G_{ad}})\otimes_{R(\tT_{comp})} PLP(\cf^+).\]
\ecor

\begin{proof}
  Since $K^0_{\tT_{comp}\times \tT_{comp}}(X^+)\simeq PLP(\cf^{+})$ as
  an $R(\tT_{comp})$-algebra by Theorem \ref{main2}, the corollary
  follows readily from Theorem \ref{relwond}. \end{proof}

\subsection{\boldmath Ordinary $K$-ring of toroidal embeddings}\label{ordtoroidal}

Note that we obtain an $R(\tT_{comp})$-algebra structure on $K^0(\bar{G_{ad}})$ is given by mapping $e^{\chi}\in R(\tT_{comp})$ to the class $[\mathcal{L}_{\chi}]\in K^0(\bar{G_{ad}})$ of the line bundle $\mathcal{L}_{\chi}$ on $\bar{G_{ad}}$.

The following theorem describes $K^0(\mathbb{X})$ as an extension of scalars of the $R(\tT_{comp})$-algebra $PLP(\cf^+)$ to the ring $K^0(\overline{G_{ad}})$. This is an extension of \cite[Theorem 4.2]{u2}.
 
\bth\label{ordkthtoroidal} We have the following isomorphism of $K^0(\bar{G_{ad}})$-algebras

\[K^0(\mathbb{X})\cong K^0(\bar{G_{ad}})\otimes_{R(\tT_{comp})} PLP(\cf^+).\]
\eeth

\begin{proof}
Since both $\mathbb{X}$ as well as $\bar{G_{ad}}$ are $\tT\times \tT$-cellular varieties, by Theorem \ref{GEQF} they are equivariantly formal for $\tG_{comp}\times \tG_{comp}$-equivariant $K$-theory. Thus we get the following isomorphisms:

\be\label{eqformaltoroidal}K^0(\mathbb{X})\cong \mathbb{Z}\otimes_{R(\tG_{comp})\otimes R(\tG_{comp})} K^0_{\tG_{comp}\times \tG_{comp}}(\mathbb{X}) \ee and

\be\label{eqformalwonderful} K^0(\bar{G_{ad}})\cong \mathbb{Z}\otimes_{R(\tG_{comp})\otimes R(\tG_{comp})} K^0_{\tG_{comp}\times \tG_{comp}}(\bar{G_{ad}}). \ee

By applying \eqref{eqformaltoroidal} and \eqref{eqformalwonderful} on the left hand side and the right hand side respectively of the isomorphism of Corollary \ref{relwondcor}, the theorem follows.
\end{proof}

%where the $R(\tT_{comp})$-algebra structure on $PLP(\cf^+)$ by Theorem \ref{main2} and %the 

\end{document}